\documentclass[a4paper, 12 pt]{article}

\usepackage[T2A]{fontenc}
\usepackage[cp1251]{inputenc}
\usepackage[english]{babel}
\usepackage[tbtags]{amsmath}
\usepackage{amsfonts,amssymb}

\usepackage{mathrsfs}

\usepackage{geometry} 
\begin{document}
\title{Numerator polynomials of the  Riordan matrices}
 \author{E. Burlachenko}
 \date{}

 \maketitle
\begin{abstract}
Riordan matrices are infinite lower triangular matrices corresponding to the certain operators in the space of formal power series. Generalized Euler polynomials ${{g}_{n}}\left( x \right)={{\left( 1-x \right)}^{n+1}}\sum\nolimits_{m=0}^{\infty }{{{p}_{n}}}\left( m \right){{x}^{m}}$, where ${{p}_{n}}\left( m \right)$ is the polynomial of degree $\le n$, are the numerator polynomials  of the generating functions of diagonals of the ordinary Riordan matrices. Generalized Narayana polynomials ${{h}_{n}}\left( x \right)={{\left( 1-x \right)}^{2n+1}}\sum\nolimits_{m=0}^{\infty }{\left( m+1 \right)...\left( m+n \right){{p}_{n}}}\left( m \right){{x}^{m}}$ are the numerator polynomials  of the generating functions of diagonals of the exponential Riordan matrices. In paper, the properties of these two types of numerator polynomials and the constructive relationships between them are considered. Separate attention is paid to the numerator polynomials of  Riordan matrices associated with the family of series $_{\left( \beta  \right)}a\left( x \right)=a\left( x{}_{\left( \beta  \right)}{{a}^{\beta }}\left( x \right) \right)$.
\end{abstract}

\section{Introduction}

Transformations in the space of formal power series over the field of real or complex numbers and the corresponding matrices are the subject of our study. We associate the rows and columns of matrices with the generating functions of their elements, i.e., formal power series. Thus, the expression $Aa\left( x \right)=b\left( x \right)$ means that the column vector multiplied by the matrix $A$ has the generating function $a\left( x \right)=\sum\nolimits_{n=0}^{\infty }{{{a}_{n}}}{{x}^{n}}$, resultant column vector has the generating function $b\left( x \right)=\sum\nolimits_{n=0}^{\infty }{{{b}_{n}}{{x}^{n}}}$.  $n$th coefficient of the series $a\left( x \right)$, $n$th row, $n$th descending diagonal and $n$th column of the matrix $A$ will be denoted  respectively by
$$\left[ {{x}^{n}} \right]a\left( x \right),    \qquad\left[ n,\to  \right]A,   \qquad[n,\searrow ]A,   \qquad A{{x}^{n}}.$$

Matrix $\left( f\left( x \right),g\left( x \right) \right)$, $n$th column of which has the generating function $f\left( x \right){{g}^{n}}\left( x \right)$, ${{g}_{0}}=0$, is called Riordan matrix (Riordan array)  [1] – [3]. It is the product of two matrices that correspond to the operators of multiplication and composition of series:
$$\left( f\left( x \right),g\left( x \right) \right)=\left( f\left( x \right),x \right)\left( 1,g\left( x \right) \right),$$
$$\left( f\left( x \right),x \right)a\left( x \right)=f\left( x \right)a\left( x \right), \qquad\left( 1,g\left( x \right) \right)a\left( x \right)=a\left( g\left( x \right) \right),$$
$$\left( f\left( x \right),g\left( x \right) \right)\left( b\left( x \right),a\left( x \right) \right)=\left( f\left( x \right)b\left( g\left( x \right) \right),a\left( g\left( x \right) \right) \right).$$
If ${{f}_{0}}\ne 0$, ${{g}_{1}}\ne 0$,  matrix $\left( f\left( x \right),g\left( x \right) \right)$ is called proper. Proper Riordan matrices form a group, called the Riordan group. 

 Matrices 
$${{\left| {{e}^{x}} \right|}^{-1}}\left( f\left( x \right),g\left( x \right) \right)\left| {{e}^{x}} \right|={{\left( f\left( x \right),g\left( x \right) \right)}_{E}},$$
where $\left| {{e}^{x}} \right|$ is the diagonal matrix: $\left| {{e}^{x}} \right|{{x}^{n}}={{{x}^{n}}}/{n!}\;$, are called exponential Riordan matrices [4], [5]. Denote $\left[ n,\to  \right]{{\left( f\left( x \right),g\left( x \right) \right)}_{E}}={{s}_{n}}\left( x \right)$. Then
$${{\left( f\left( x \right),g\left( x \right) \right)}_{E}}{{\left( 1-\varphi x \right)}^{-1}}={{\left| {{e}^{x}} \right|}^{-1}}\left( f\left( x \right),g\left( x \right) \right){{e}^{\varphi x}}={{\left| {{e}^{x}} \right|}^{-1}}f\left( x \right)\exp \left( \varphi g\left( x \right) \right),$$
or
$$\sum\limits_{n=0}^{\infty }{\frac{{{s}_{n}}\left( \varphi  \right)}{n!}{{x}^{n}}}=f\left( x \right)\exp \left( \varphi g\left( x \right) \right).$$
If ${{f}_{0}}\ne 0$, ${{g}_{1}}\ne 0$, sequence of polynomials ${{s}_{n}}\left( x \right)$ is called Sheffer sequence, and in the case $f\left( x \right)=1$ binomial sequence. Properties of the Sheffer sequences are subject of study of the umbral calculus [6]. Matrix
$$P=\left( \frac{1}{1-x},\frac{x}{1-x} \right)={{\left( {{e}^{x}},x \right)}_{E}}=\left( \begin{matrix}
   1 & 0 & 0 & 0 & \cdots   \\
   1 & 1 & 0 & 0 & \cdots   \\
   1 & 2 & 1 & 0 & \cdots   \\
   1 & 3 & 3 & 1 & \cdots   \\
   \vdots  & \vdots  & \vdots  & \vdots  & \ddots   \\
\end{matrix} \right)$$
is called Pascal matrix. Power of the Pascal matrix is defined by 
$${{P}^{\varphi }}=\left( \frac{1}{1-\varphi x},\frac{x}{1-\varphi x} \right)={{\left( {{e}^{\varphi x}},x \right)}_{E}}.$$

Along with the lower triangular Riordan matrices, we will consider the “square” matrices $\left( b\left( x \right),a\left( x \right) \right)$,  ${{b}_{0}}\ne 0$, ${{a}_{0}}=1$, whose $n$th column has the generating function $b\left( x \right){{a}^{n}}\left( x \right)$. Upper triangular matrix $\left( 1,1+x \right)$, transposed to the Pascal matrix and coinciding with the matrix of shift operator, also belongs to this type of matrices:
$$\left( 1,1+x \right)={{P}^{T}}=E=\left( \begin{matrix}
   1 & 1 & 1 & 1 & \cdots   \\
   0 & 1 & 2 & 3 & \cdots   \\
   0 & 0 & 1 & 3 & \cdots   \\
   0 & 0 & 0 & 1 & \cdots   \\
   \vdots  & \vdots  & \vdots  & \vdots  & \ddots   \\
\end{matrix} \right).$$
Matrix $\left( b\left( x \right),a\left( x \right) \right)$ can be multiplied from the right by the matrix with the finite columns and from the left by the matrix with the finite rows. Since $\left[ n,\to  \right]\left( b\left( x \right),a\left( x \right) \right)=\linebreak=[n,\searrow ]\left( b\left( x \right),xa\left( x \right) \right)$, then the matrix$\left( b\left( x \right),a\left( x \right) \right)$ is a useful tool for studying the matrix $\left( b\left( x \right),xa\left( x \right) \right)$.  Denote
$$\left[ n,\to  \right]\left( b\left( x \right),a\left( x \right)-1 \right)={{w}_{n}}\left( x \right)=\sum\limits_{m=0}^{n}{{{w}_{m}}{{x}^{m}}}.$$
Since
$$\left( b\left( x \right),a\left( x \right) \right)=\left( b\left( x \right),a\left( x \right)-1 \right)\left( 1,1+x \right),  \qquad\left[ n,\to  \right]\left( 1,1+x \right)=\frac{{{x}^{n}}}{{{\left( 1-x \right)}^{n+1}}},$$
then
$$\left[ n,\to  \right]\left( b\left( x \right),a\left( x \right) \right)=\sum\limits_{m=0}^{n}{\frac{{{w}_{m}}{{x}^{m}}}{{{\left( 1-x \right)}^{m+1}}}}=\frac{\sum\limits_{m=0}^{n}{{{w}_{m}}{{x}^{m}}{{\left( 1-x \right)}^{n-m}}}}{{{\left( 1-x \right)}^{n+1}}}=\frac{{{g}_{n}}\left( x \right)}{{{\left( 1-x \right)}^{n+1}}}.$$
If $b\left( x \right)=1$, $a\left( x \right)={{e}^{x}}$, then ${{g}_{n}}\left( x \right)={{{A}_{n}}\left( x \right)}/{n!}\;$, where ${{A}_{n}}\left( x \right)$ are the Euler polynomials:
$$\frac{{{A}_{n}}\left( x \right)}{{{\left( 1-x \right)}^{n+1}}}=\sum\limits_{m=0}^{\infty }{{{m}^{n}}}{{x}^{m}},  \qquad{{A}_{n}}\left( 1 \right)=n!.$$
For example,
$${{A}_{0}}\left( x \right)=1,  \qquad{{A}_{1}}\left( x \right)=x, \qquad{{A}_{2}}\left( x \right)=x+{{x}^{2}},  \qquad{{A}_{3}}\left( x \right)=x+4{{x}^{2}}+{{x}^{3}},$$
$${{A}_{4}}\left( x \right)=x+11{{x}^{2}}+11{{x}^{3}}+{{x}^{4}}.$$
In this connection, we will called polynomials ${{g}_{n}}\left( x \right)$ the generalized Euler polynomials. 

“Square” Riordan matrices (called convolution arrays) and numerator polynomials of the generating functions of their rows were considered in the series of papers [7] – [11]. In [12] such matrices are called generalized Riordan arrays. Concept of generalized Euler polynomials (called ${{p}_{n}}$-associated Eulerian polynomials) in general form   is represented in [13]. 

Mentioned series of papers was written before the advent of the Riordan matrices theory. In modern literature, the topic of  numerator polynomials of the  Riordan matrices has not been raised. Exception consists of the three works [14] – [16], in which the generating functions of the diagonals of exponential Riordan matrices are considered. These works provided an incentive for writing the present paper.

In paper, we study the transformations in which the numerator polynomials of  Riordan matrices participate, and build matrices of these transformations according to certain rules. These matrices serve as tools for further research. In Section 2, we consider the numerator polynomials of ordinary Riordan matrices (generalized Euler polynomials) and the transformations associated with them. Using matrices of these transformations, we obtain a general formula for the numerator polynomials of the matrices $\left( 1,{{x}_{\left( \beta  \right)}}a\left( x \right) \right)$, where $_{\left( \beta  \right)}a\left( x \right)$ is the generalized binomial series:
$$_{\left( \beta  \right)}{{a}^{\varphi }}\left( x \right)=\sum\limits_{n=0}^{\infty }{\frac{\varphi }{\varphi +\beta n}}\left( \begin{matrix}
   \varphi +\beta n  \\
   n  \\
\end{matrix} \right){{x}^{n}}.$$
In Section 3, we consider the numerator polynomials of exponential Riordan matrices (generalized Narayana polynomials) and the transformations associated with them. In Section 4, we introduce matrices of the operators mapping generalized Euler polynomials to generalized Narayana polynomials. Using these matrices, we obtain a general formula for the numerator polynomials of the matrices ${{\left( 1,{{x}_{\left( \beta  \right)}}a\left( x \right) \right)}_{E}}$, where $_{\left( \beta  \right)}a\left( x \right)$ is the generalized binomial series. In Section 5, we give an idea of the generalized Lagrange series $_{\left( \beta  \right)}a\left( x \right)$ defined by the identity $_{\left( \beta  \right)}a\left( x \right)=a\left( x{}_{\left( \beta  \right)}{{a}^{\beta }}\left( x \right) \right)$, or ${{\left( 1,x{{a}^{-\beta }}\left( x \right) \right)}^{-1}}=\left( 1,x{}_{\left( \beta  \right)}{{a}^{\beta }}\left( x \right) \right)$. In other words, we extend the generalization that underlies the generalized binomial series to any formal power series $a\left( x \right)$, ${{a}_{0}}=1$. In the remainder of paper, we will consider the transformations associated with these series.
\section{ Generalized Euler polynomials }
 In [13], the transformations associated with the generalized Euler polynomials are presented from a general point of view. In this section, we consider them from  point of view of the Riordan matrix theory. We denote
$$\left[ n,\to  \right]\left( b\left( x \right),\log a\left( x \right) \right)={{c}_{n}}\left( x \right)=\sum\limits_{m=0}^{n}{{{c}_{m}}{{x}^{m}}}, \qquad{{b}_{0}}\ne 0, \qquad{{a}_{0}}=1,$$
$$\left[ n,\to  \right]{{\left( b\left( x \right),\log a\left( x \right) \right)}_{E}}={{s}_{n}}\left( x \right)=\sum\limits_{m=0}^{n}{{{s}_{m}}{{x}^{m}}},
\qquad\left[ n,\to  \right]\left( b\left( x \right),a\left( x \right) \right)=\frac{{{g}_{n}}\left( x \right)}{{{\left( 1-x \right)}^{n+1}}}.$$
Then
$$b\left( x \right){{a}^{m}}\left( x \right)=\sum\limits_{n=0}^{\infty }{\frac{{{s}_{n}}\left( m \right)}{n!}{{x}^{n}}},  \qquad\frac{{{g}_{n}}\left( x \right)}{{{\left( 1-x \right)}^{n+1}}}=\sum\limits_{m=0}^{\infty }{\frac{{{s}_{n}}\left( m \right)}{n!}{{x}^{m}}}.$$  
Since
$$\left( b\left( x \right),a\left( x \right) \right)=\left( b\left( x \right),\log a\left( x \right) \right)\left( 1,{{e}^{x}} \right),$$
then
$$\frac{{{g}_{n}}\left( x \right)}{{{\left( 1-x \right)}^{n+1}}}=\sum\limits_{p=0}^{n}{\frac{{{{c}_{p}}{{A}_{p}}\left( x \right)}/{p!}\;}{{{\left( 1-x \right)}^{p+1}}}}=\frac{1}{n!}\sum\limits_{p=0}^{n}{\frac{{{s}_{p}}{{A}_{p}}\left( x \right)}{{{\left( 1-x \right)}^{p+1}}}}=\frac{\frac{1}{n!}\sum\limits_{p=0}^{n}{{{s}_{p}}{{\left( 1-x \right)}^{n-p}}{{A}_{p}}\left( x \right)}}{{{\left( 1-x \right)}^{n+1}}}.$$
We introduce the matrices ${{U}_{n}}$:
$${{U}_{n}}{{x}^{p}}={{\left( 1-x \right)}^{n+1}}\frac{1}{n!}\sum\limits_{m=0}^{\infty }{{{m}^{p}}}{{x}^{m}}=\frac{1}{n!}{{\left( 1-x \right)}^{n-p}}{{A}_{p}}\left( x \right),  \qquad 0\le p\le n.$$
For example,
$${{U}_{1}}=\left( \begin{matrix}
   1 & 0  \\
   -1 & 1  \\
\end{matrix} \right), \qquad{{U}_{2}}=\frac{1}{2!}\left( \begin{matrix}
   1 & 0 & 0  \\
   -2 & 1 & 1  \\
   1 & -1 & 1  \\
\end{matrix} \right), \qquad{{U}_{3}}=\frac{1}{3!}\left( \begin{matrix}
   1 & 0 & 0 & 0  \\
   -3 & 1 & 1 & 1  \\
   3 & -2 & 0 & 4  \\
   -1 & 1 & -1 & 1  \\
\end{matrix} \right).$$
Then ${{U}_{n}}{{s}_{n}}\left( x \right)={{g}_{n}}\left( x \right)$. Let the symbols ${{\left( \varphi  \right)}_{n}}$, ${{\left[ \varphi  \right]}_{n}}$ denote respectively the descending and ascending  factorial:
$${{\left( \varphi  \right)}_{n}}=\varphi \left( \varphi -1 \right)...\left( \varphi -n+1 \right),   \qquad{{\left[ \varphi  \right]}_{n}}=\varphi \left( \varphi +1 \right)...\left( \varphi +n-1 \right).$$
Since
$$\frac{{{x}^{p}}}{{{\left( 1-x \right)}^{n+1}}}=\sum\limits_{m=0}^{\infty }{\left( \begin{matrix}
   m+n-p  \\
   n  \\
\end{matrix} \right)}{{x}^{m}}=\sum\limits_{m=0}^{\infty }{\frac{{{\left( m+n-p \right)}_{n}}}{n!}}{{x}^{m}}, \qquad 0\le p\le n,$$
then
$$U_{n}^{-1}{{x}^{p}}={{\left( x+n-p \right)}_{n}}={{\left( x \right)}_{p}}{{\left[ x+1 \right]}_{n-p}}.$$
For example,
$$U_{1}^{-1}=\left( \begin{matrix}
   1 & 0  \\
   1 & 1  \\
\end{matrix} \right),  \qquad U_{2}^{-1}=\left( \begin{matrix}
   2 & 0 & 0  \\
   3 & 1 & -1  \\
   1 & 1 & 1  \\
\end{matrix} \right), \qquad U_{3}^{-1}=\left( \begin{matrix}
   6 & 0 & 0 & 0  \\
   11 & 2 & -1 & 2  \\
   6 & 3 & 0 & -3  \\
   1 & 1 & 1 & 1  \\
\end{matrix} \right).$$

We will consider the products of matrices of various orders with infinite matrices. Thus, the matrix of order $n$ corresponds to infinite matrix, the columns elements of which, starting from the $n$th column, $n=0$, $1$,$2$, $...$, and the rows elements of which, starting from the $n$th row, are equal to zero. We introduce the matrices ${{J}_{n}}$ corresponding to the operator rearranging a coefficients of  polynomial of degree $n$ in the reverse order: ${{J}_{n}}c\left( x \right)={{x}^{n}}c\left( {1}/{x}\; \right)$, where $c\left( x \right)$ is a polynomial of degree $\le n$. For example,
$${{J}_{3}}=\left( \begin{matrix}
   0 & 0 & 0 & 1  \\
   0 & 0 & 1 & 0  \\
   0 & 1 & 0 & 0  \\
   1 & 0 & 0 & 0  \\
\end{matrix} \right).$$
{\bfseries Theorem 2.1.} \emph{
$${{U}_{n}}E\left( 1,-x \right)U_{n}^{-1}={{\left( -1 \right)}^{n}}{{J}_{n}}.$$}
{\bfseries Proof.}
$$E\left( 1,-x \right){{\left( x \right)}_{p}}{{\left[ x+1 \right]}_{n-p}}={{\left( -x-1 \right)}_{p}}{{\left[ -x \right]}_{n-p}}={{\left( -1 \right)}^{n}}{{\left( x \right)}_{n-p}}{{\left[ x+1 \right]}_{p}},$$
or
$$E\left( 1,-x \right)U_{n}^{-1}={{\left( -1 \right)}^{n}}U_{n}^{-1}{{J}_{n}}. \qquad\square $$ 
{\bfseries Theorem 2.2.} \emph{Polynomials ${{\left( -1 \right)}^{n}}{{J}_{n}}{{g}_{n}}\left( x \right)$ are the numerator polynomials of the  matrix $\left( b\left( x \right){{a}^{-1}}\left( x \right),x{{a}^{-1}}\left( x \right) \right)$.}
\\{\bfseries Proof.} By the Theorem 2.1.
$${{\left( -1 \right)}^{n}}{{J}_{n}}{{g}_{n}}\left( x \right)={{U}_{n}}E\left( 1,-x \right){{s}_{n}}\left( x \right).$$ 
Since
$$\left( b\left( x \right),\log a\left( x \right) \right)\left( 1,-x \right)\left( {{e}^{x}},x \right)=\left( b\left( x \right){{a}^{-1}}\left( x \right),\log {{a}^{-1}}\left( x \right) \right),$$
then
$$E\left( 1,-x \right)s\left( x \right)={{s}_{n}}\left( -x-1 \right)=\left[ n,\to  \right]{{\left( b\left( x \right){{a}^{-1}}\left( x \right),\log {{a}^{-1}}\left( x \right) \right)}_{E}}. \qquad\square $$
{\bfseries Theorem 2.3.} \emph{$${{g}_{n}}\left( 1 \right)={{b}_{0}}{{\left( {{a}_{1}} \right)}^{n}}.$$}
{\bfseries Proof.} Denote ${{U}_{n}}{{x}^{p}}={{U}_{n,p}}\left( x \right)$. Since
$${{a}_{1}}=\left[ x \right]\log a\left( x \right), \qquad{{b}_{0}}{{\left( {{a}_{1}} \right)}^{n}}=\left[ {{x}^{n}} \right]{{s}_{n}}\left( x \right);  \qquad{{U}_{n,p}}\left( 1 \right)=0, \qquad p<n;  \qquad{{U}_{n,n}}\left( 1 \right)=1,$$
then
$${{g}_{n}}\left( 1 \right)=\sum\limits_{p=0}^{n}{{{s}_{p}}{{U}_{n,p}}}\left( 1 \right)={{s}_{n}}={{b}_{0}}{{\left( {{a}_{1}} \right)}^{n}}. \qquad\square $$

Case, when ${{a}_{1}}=0$ and  the degree of  polynomial ${{s}_{n}}\left( x \right)$ is less than $n$, is possible. This possibility is reflected in the following theorem.
\\{\bfseries Theorem 2.4.} \emph{If ${{c}_{n-m}}\left( x \right)$ is a polynomial of degree $n-m$, then
$${{U}_{n}}{{c}_{n-m}}\left( x \right)={{\left( 1-x \right)}^{m}}\frac{\left( n-m \right)!}{n!}{{U}_{n-m}}{{c}_{n-m}}\left( x \right).$$
Respectively, if ${{d}_{n-m}}\left( x \right)$ is a polynomial of degree $\le n-m$, then
$$U_{n}^{-1}{{\left( 1-x \right)}^{m}}{{d}_{n-m}}\left( x \right)=\frac{n!}{\left( n-m \right)!}U_{n-m}^{-1}{{d}_{n-m}}\left( x \right).$$}
{\bfseries Proof.} Let ${{I}_{n}}$ is the identity matrix of order $n+1$. It is obvious that
$$\left( {{\left( 1-x \right)}^{-m}},x \right){{U}_{n}}{{I}_{n-m}}=\frac{\left( n-m \right)!}{n!}{{U}_{n-m}},$$
or
$${{U}_{n}}{{I}_{n-m}}=\left( {{\left( 1-x \right)}^{m}},x \right)\frac{\left( n-m \right)!}{n!}{{U}_{n-m}}.$$
Respectively,
 $$U_{n}^{-1}\left( {{\left( 1-x \right)}^{m}},x \right){{I}_{n-m}}=\frac{n!}{\left( n-m \right)!}U_{n-m}^{-1}.\qquad \square $$

We introduce  the matrices  ${{V}_{n}}={{J}_{n}}E{{J}_{n}}=\left( {{\left( 1+x \right)}^{n+1}},x \right){{P}^{-1}}{{I}_{n}}$. For example,
$${{V}_{3}}=\left( \begin{matrix}
   1 & 0 & 0 & 0  \\
   3 & 1 & 0 & 0  \\
   3 & 2 & 1 & 0  \\
   1 & 1 & 1 & 1  \\
\end{matrix} \right),  \qquad V_{3}^{-1}=\left( \begin{matrix}
   1 & 0 & 0 & 0  \\
   -3 & 1 & 0 & 0  \\
   3 & -2 & 1 & 0  \\
   -1 & 1 & -1 & 1  \\
\end{matrix} \right).$$
If $c\left( x \right)$ is a polynomial of degree $\le n$, then
$${{V}_{n}}c\left( x \right)={{\left( 1+x \right)}^{n}}c\left( \frac{x}{1+x} \right),  \qquad V_{n}^{-1}c\left( x \right)={{\left( 1-x \right)}^{n}}c\left( \frac{x}{1-x} \right).$$
We found in the introduction that $V_{n}^{-1}{{w}_{n}}\left( x \right)={{g}_{n}}\left( x \right)$, where ${{w}_{n}}\left( x \right)=\left[ n,\to  \right]\left( b\left( x \right),a\left( x \right)-1 \right)$, and hence $U_{n}^{-1}V_{n}^{-1}{{w}_{n}}\left( x \right)={{s}_{n}}\left( x \right)$. By the Theorem 2.4. we find:
$$U_{n}^{-1}V_{n}^{-1}{{x}^{p}}=U_{n}^{-1}{{\left( 1-x \right)}^{n-p}}{{x}^{p}}=\frac{n!}{p!}{{\left( x \right)}_{p}}=\frac{n!}{p!}\sum\limits_{m=0}^{p}{s\left( p,\text{ }m \right){{x}^{m}}},$$
where $s\left( p,\text{ }m \right)$ are the Stirling numbers of  the first kind. Hence
$${{V}_{n}}{{U}_{n}}{{x}^{p}}=\frac{1}{n!}\sum\limits_{m=0}^{p}{m!S\left( p,\text{ }m \right)}\text{ }{{x}^{m}},$$
where $S\left( p,\text{ }m \right)$ are the Stirling numbers of the second kind. For example, $U_{3}^{-1}V_{3}^{-1}$, ${{V}_{3}}{{U}_{3}}$:
$$3!\left( \begin{matrix}
   1 & 0 & 0 & 0  \\
   0 & 1 & -1 & 2  \\
   0 & 0 & 1 & -3  \\
   0 & 0 & 0 & 1  \\
\end{matrix} \right)\left( \begin{matrix}
   1 & 0 & 0 & 0  \\
   0 & 1 & 0 & 0  \\
   0 & 0 & \frac{1}{2!} & 0  \\
   0 & 0 & 0 & \frac{1}{3!}  \\
\end{matrix} \right),  \qquad\frac{1}{3!}\left( \begin{matrix}
   1 & 0 & 0 & 0  \\
   0 & 1 & 0 & 0  \\
   0 & 0 & 2! & 0  \\
   0 & 0 & 0 & 3!  \\
\end{matrix} \right)\left( \begin{matrix}
   1 & 0 & 0 & 0  \\
   0 & 1 & 1 & 1  \\
   0 & 0 & 1 & 3  \\
   0 & 0 & 0 & 1  \\
\end{matrix} \right).$$

Numerator polynomials of the matrix $\left( 1,xa\left( x \right) \right)$ are of independent interest. We introduce special notation for them and the associated polynomials:
$$\left[ n,\searrow  \right]\left( 1,xa\left( x \right) \right)=\frac{{{\alpha }_{n}}\left( x \right)}{{{\left( 1-x \right)}^{n+1}}},  \qquad\left[ n,\searrow  \right]\left( 1,x{{a}^{-1}}\left( x \right) \right)=\frac{\alpha _{n}^{\left( -1 \right)}\left( x \right)}{{{\left( 1-x \right)}^{n+1}}},$$
$$\left[ n,\to  \right]{{\left( 1,\log a\left( x \right) \right)}_{E}}={{u}_{n}}\left( x \right),  \qquad\left[ n,\to  \right]\left( 1,a\left( x \right)-1 \right)={{v}_{n}}\left( x \right).$$
Since $\left[ n,\searrow  \right]\left( {{a}^{-1}}\left( x \right),x{{a}^{-1}}\left( x \right) \right)=\left( {1}/{x}\; \right)\left[ n,\searrow  \right]\left( 1,x{{a}^{-1}}\left( x \right) \right)$, $n>0$, then by  the \linebreak Theorem 2. 2.
 $$\alpha _{n}^{\left( -1 \right)}\left( x \right)={{\left( -1 \right)}^{n}}x{{J}_{n}}{{\alpha }_{n}}\left( x \right).$$
{\bfseries Example 2.1.}
$$a\left( x \right)=\frac{1+x}{1-x},   \qquad a\left( x \right)-1=\frac{2x}{1-x},  \qquad{{v}_{n}}\left( x \right)={{2}^{n}}x{{\left( \frac{1}{2}+x \right)}^{n-1}}, \qquad n>0,$$
$${{\alpha }_{n}}\left( x \right)=V_{n}^{-1}{{v}_{n}}\left( x \right)=2x{{\left( 1+x \right)}^{n-1}},$$
$${{u}_{n}}\left( x \right)=U_{n}^{-1}{{\alpha }_{n}}\left( x \right)=2\sum\limits_{p=0}^{n}{\left( \begin{matrix}
   n-1  \\
   p-1  \\
\end{matrix} \right)}{{\left( x \right)}_{p}}{{\left[ x+1 \right]}_{n-p}}=$$
$$=U_{n}^{-1}V_{n}^{-1}{{v}_{n}}\left( x \right)=n!\sum\limits_{p=0}^{n}{\left( \begin{matrix}
   n-1  \\
   p-1  \\
\end{matrix} \right)}\frac{{{2}^{p}}}{p!}{{\left( x \right)}_{p}}.$$
{\bfseries Example 2.2.} Let ${{\left( 1,xg\left( x \right) \right)}^{-1}}=\left( 1,x{{b}^{-1}}\left( x \right) \right)$, ${{g}_{0}}=1$, ${{b}_{0}}=1$. Then by the Lagrange inversion theorem
$$\left[ {{x}^{n}} \right]{{g}^{m}}\left( x \right)=\frac{m}{m+n}\left[ {{x}^{n}} \right]{{b}^{m+n}}\left( x \right)=\left[ {{x}^{n}} \right]\left( 1-x{{\left( \log b\left( x \right) \right)}^{\prime }} \right){{b}^{m+n}}\left( x \right).$$
Denote $\left[ {{x}^{n}} \right]{{g}^{m}}\left( x \right)=g_{n}^{\left( m \right)}$, $\left[ {{x}^{n}} \right]\left( 1-x{{\left( \log b\left( x \right) \right)}^{\prime }} \right){{b}^{m}}\left( x \right)=c_{n}^{\left( m \right)}$. Then
$$\left[ n,\to  \right]\left( 1,xg\left( x \right) \right)=\sum\limits_{m=0}^{n}{g_{n-m}^{\left( m \right)}}{{x}^{m}}=\sum\limits_{m=0}^{n}{c_{n-m}^{\left( n \right)}{{x}^{m}}}=\left[ n,\to  \right]\left( \left( 1-x{{\left( \log b\left( x \right) \right)}^{\prime }} \right){{b}^{n}}\left( x \right),x \right).$$
If
$$g\left( x \right)=\frac{x}{2}+\sqrt{1+\frac{{{x}^{2}}}{4}},$$
then
$$b\left( x \right)=\sqrt{1+x},    \qquad1-x{{\left( \log b\left( x \right) \right)}^{\prime }}=\left( 1+\frac{x}{2} \right)\frac{1}{1+x},$$
$$\left[ 2n,\to  \right]\left( 1,xg\left( x \right) \right)=\left( \frac{1}{2}+x \right){{x}^{n}}{{\left( 1+x \right)}^{n-1}},  \qquad n>0.$$
Let
$$a\left( x \right)={{\left( \frac{x}{2}+\sqrt{1+\frac{{{x}^{2}}}{4}} \right)}^{2}}.$$
Then
$$a\left( x \right)-1=x\sqrt{a\left( x \right)}=x\left( \frac{x}{2}+\sqrt{1+\frac{{{x}^{2}}}{4}} \right),$$  
$${{v}_{2n}}\left( x \right)=\left( \frac{1}{2}+x \right){{x}^{n}}{{\left( 1+x \right)}^{n-1}},  \qquad{{\alpha }_{2n}}\left( x \right)=V_{2n}^{-1}{{v}_{2n}}\left( x \right)=\frac{1}{2}\left( 1+x \right){{x}^{n}},$$
$${{u}_{2n}}\left( x \right)=U_{2n}^{-1}{{\alpha }_{2n}}\left( x \right)=\frac{1}{2}{{\left( x \right)}_{n}}{{\left[ x+1 \right]}_{n}}+\frac{1}{2}{{\left( x \right)}_{n+1}}{{\left[ x+1 \right]}_{n-1}}=\prod\limits_{m=0}^{n-1}{\left( {{x}^{2}}-{{m}^{2}} \right)}.$$

We find the generating function of the sequence of polynomials ${{\alpha }_{n}}\left( t \right)$. Since
$$\frac{{{\alpha }_{n}}\left( t \right)}{{{\left( 1-t \right)}^{n+1}}}=\sum\limits_{m=0}^{\infty }{\left[ {{x}^{n}} \right]{{a}^{m}}\left( x \right){{t}^{m}}}=\left[ {{x}^{n}} \right]\frac{1}{1-ta\left( x \right)},$$
then 
$$\sum\limits_{n=0}^{\infty }{{{\alpha }_{n}}}\left( t \right){{x}^{n}}=\frac{1-t}{1-ta\left( x\left( 1-t \right) \right)}.\eqno(1)$$	
{\bfseries Example 2.3.} This example is a generalization of the examples considered in [11]. Denote ${{\tilde{\alpha }}_{0}}\left( x \right)=1$, ${{\tilde{\alpha }}_{n}}\left( x \right)=\left( {1}/{x}\; \right){{\alpha }_{n}}\left( x \right)$.  If
$$\left[ n,\searrow  \right]\left( 1,\frac{x}{1+\varphi x+\beta {{x}^{2}}} \right)=\frac{{{\alpha }_{n}}\left( x \right)}{{{\left( 1-x \right)}^{n+1}}},$$ 
then
$${{\tilde{\alpha }}_{n}}\left( x \right)=\left[ n,\to  \right]\left( \frac{1}{1+\varphi x+\beta {{x}^{2}}},\frac{\beta {{x}^{2}}}{1+\varphi x+\beta {{x}^{2}}} \right).$$
Really,
$$1+\sum\limits_{n=1}^{\infty }{{{{\tilde{\alpha }}}_{n}}}\left( t \right){{x}^{n}}=\left( \frac{1}{1+\varphi x+\beta {{x}^{2}}},\frac{\beta {{x}^{2}}}{1+\varphi x+\beta {{x}^{2}}} \right)\frac{1}{1-tx}=
1+\frac{-\varphi x-\beta \left( 1-t \right){{x}^{2}}}{1+\varphi x+\beta \left( 1-t \right){{x}^{2}}},$$
$$\sum\limits_{n=0}^{\infty }{{{\alpha }_{n}}}\left( t \right){{x}^{n}}=1+t\sum\limits_{n=1}^{\infty }{{{{\tilde{\alpha }}}_{n}}}\left( t \right){{x}^{n}}=\frac{1+\varphi \left( 1-t \right)x+\beta {{\left( 1-t \right)}^{2}}{{x}^{2}}}{1+\varphi x+\beta \left( 1-t \right){{x}^{2}}},$$ 
that corresponds to the formula (1).

Let $_{\left( \beta  \right)}a\left( x \right)$ is the generalized binomial series: 
$$_{\left( \beta  \right)}{{a}^{\varphi }}\left( x \right)=\sum\limits_{n=0}^{\infty }{\frac{\varphi }{\varphi +n\beta }}\left( \begin{matrix}
   \varphi +n\beta   \\
   n  \\
\end{matrix} \right){{x}^{n}};$$
$$_{\left( 0 \right)}a\left( x \right)=1+x,  \qquad_{\left( 1 \right)}a\left( x \right)=\frac{1}{1-x},    \qquad_{\left( 2 \right)}a\left( x \right)=\frac{1-\sqrt{1-4x}}{2x},$$
$$_{\left( -1 \right)}a\left( x \right)=\frac{1+\sqrt{1+4x}}{2},    \qquad_{\left( {1}/{2}\; \right)}a\left( x \right)={{\left( \frac{x}{2}+\sqrt{1+\frac{{{x}^{2}}}{4}} \right)}^{2}}.$$
Denote
$$\left[ n,\searrow  \right]\left( 1,{{x}_{\left( \beta  \right)}}a\left( x \right) \right)=\frac{_{\left( \beta  \right)}{{\alpha }_{n}}\left( x \right)}{{{\left( 1-x \right)}^{n+1}}}=\sum\limits_{m=0}^{\infty }{\frac{m}{m+n\beta }}\left( \begin{matrix}
   m+n\beta   \\
   n  \\
\end{matrix} \right){{x}^{m}},$$
$$\left[ n,\to  \right]\left( {{1,}_{\left( \beta  \right)}}a\left( x \right)-1 \right)={}_{\left( \beta  \right)}{{v}_{n}}\left( x \right).$$
{\bfseries Theorem 2.5.} \emph{}
$$_{\left( \beta  \right)}{{\alpha }_{n}}\left( x \right)=\frac{1}{n}\sum\limits_{m=0}^{n}{\left( \begin{matrix}
   n\left( 1-\beta  \right)  \\
   m-1  \\
\end{matrix} \right)\left( \begin{matrix}
   n\beta   \\
   n-m  \\
\end{matrix} \right){{x}^{m}}}, \qquad n>0.\eqno(2)$$
{\bfseries Proof.} We take into account the property of generalized binomial series: $_{\left( \beta  \right)}a\left( x \right)-1={{x}_{\left( \beta  \right)}}{{a}^{\beta }}\left( x \right)$. Then
$${}_{\left( \beta  \right)}{{v}_{n}}\left( x \right)=\left[ n,\to  \right]\left( 1,{{x}_{\left( \beta  \right)}}{{a}^{\beta }}\left( x \right) \right)=\sum\limits_{m=0}^{n}{\frac{m}{n}\left( \begin{matrix}
   n\beta   \\
   n-m  \\
\end{matrix} \right)}{{x}^{m}}.$$
We use the factorial representation of binomial coefficients, i.e., we prove the theorem for positive integers $\beta $. By polynomial argument (binomial coefficients under consideration are polynomials in $\beta $) this is equivalent to the general proof. Since 
$$\left[ m,\to  \right]V_{n}^{-1}=\sum\limits_{i=0}^{m}{\left( \begin{matrix}
   m-n-1  \\
   m-i  \\
\end{matrix} \right){{x}^{i}}}=\sum\limits_{i=0}^{m}{{{\left( -1 \right)}^{m-i}}}\left( \begin{matrix}
   n-i  \\
   m-i  \\
\end{matrix} \right){{x}^{i}},$$
then 
$$\left[ {{x}^{m}} \right]{}_{\left( \beta  \right)}{{\alpha }_{n}}\left( x \right)=\left[ {{x}^{m}} \right]V{{_{n}^{-1}}_{\left( \beta  \right)}}{{v}_{n}}\left( x \right)=$$
$$=\sum\limits_{i=0}^{m}{{{\left( -1 \right)}^{m-i}}}\left( \begin{matrix}
   n-i  \\
   m-i  \\
\end{matrix} \right)\frac{i}{n}\left( \begin{matrix}
   n\beta   \\
   n-i  \\
\end{matrix} \right)\frac{\left( n\beta -n+m \right)!}{\left( n\beta -n+m \right)!}=$$
$$=\frac{1}{n}\left( \begin{matrix}
   n\beta   \\
   n-m  \\
\end{matrix} \right)\sum\limits_{i=0}^{m}{{{\left( -1 \right)}^{m-i}}i}\left( \begin{matrix}
   n\beta -n+m  \\
   m-i  \\
\end{matrix} \right)=$$ 
$$=\frac{1}{n}\left( \begin{matrix}
   n\beta   \\
   n-m  \\
\end{matrix} \right){{\left( -1 \right)}^{m-1}}\left( \begin{matrix}
   n\beta -n+m-2  \\
   m-1  \\
\end{matrix} \right)=\frac{1}{n}\left( \begin{matrix}
   n\beta   \\
   n-m  \\
\end{matrix} \right)\left( \begin{matrix}
   n\left( 1-\beta  \right)  \\
   m-1  \\
\end{matrix} \right).\qquad\square $$ 

Note that
$$_{\left( 0 \right)}{{\alpha }_{n}}\left( x \right)={{x}^{n}}, \qquad_{\left( 1 \right)}{{\alpha }_{n}}\left( x \right)=x,  \qquad_{\left( {1}/{2}\; \right)}{{\alpha }_{2n}}\left( x \right)=\frac{1}{2}\left( 1+x \right){{x}^{n}}.$$
Since $_{\left( 1-\beta  \right)}a\left( x \right)={}_{\left( \beta  \right)}{{a}^{-1}}\left( -x \right)$, then $_{\left( 1-\beta  \right)}{{\alpha }_{n}}\left( x \right)=x{{J}_{n}}{}_{\left( \beta  \right)}{{\alpha }_{n}}\left( x \right)$.
\section{Generalized Narayana polynomials }

Constructive relationships between the ordinary and exponential Riordan matrices exist. Particular manifestations of these relationships resemble the details of construction, the general plan of which is a secret for us. Following [14] – [16], we will consider some of such manifestations associated with the numerator polynomials. Since
$$[n,\searrow ]\left( b\left( x \right),xa\left( x \right) \right)=\sum\limits_{m=0}^{\infty }{\frac{{{s}_{n}}\left( m \right)}{n!}}{{x}^{m}},$$ 
then
$$[n,\searrow ]{{\left( b\left( x \right),xa\left( x \right) \right)}_{E}}=\sum\limits_{m=0}^{\infty }{\frac{\left( m+n \right)!}{m!}\frac{{{s}_{n}}\left( m \right)}{n!}}{{x}^{m}}=\sum\limits_{m=0}^{\infty }{\frac{{{\left[ m+1 \right]}_{n}}{{s}_{n}}\left( m \right)}{n!}{{x}^{m}}}.$$  
  If ${{a}_{1}}\ne 0$,  then ${{\left[ x+1 \right]}_{n}}{{s}_{n}}\left( x \right)$ is the polynomial of degree $2n$, in general case –  of degree $\le 2n$, so that
$$\sum\limits_{m=0}^{\infty }{\frac{{{\left[ m+1 \right]}_{n}}{{s}_{n}}\left( m \right)}{n!}{{x}^{m}}}=\frac{{{h}_{n}}\left( x \right)}{{{\left( 1-x \right)}^{2n+1}}},   \qquad{{h}_{n}}\left( x \right)=\frac{\left( 2n \right)!}{n!}{{U}_{2n}}{{\left[ x+1 \right]}_{n}}{{s}_{n}}\left( x \right).$$
Since ${{\left[ -m \right]}_{n}}=0$ when $m=0$, $1$, … , $n-1$, then in accordance with the Theorem 2.2.
$$\sum\limits_{m=0}^{\infty }{\frac{{{\left[ -m \right]}_{n}}{{s}_{n}}\left( -m-1 \right)}{n!}{{x}^{m}}}=\frac{{{\left( -1 \right)}^{2n}}{{J}_{2n}}{{h}_{n}}\left( x \right)}{{{\left( 1-x \right)}^{2n+1}}}=\frac{{{\left( -1 \right)}^{2n}}{{x}^{n}}{{J}_{n}}{{h}_{n}}\left( x \right)}{{{\left( 1-x \right)}^{2n+1}}},$$
i.e., ${{h}_{n}}\left( x \right)$  is the polynomial of degree $\le n$. Since $\left[ {{x}^{2n}} \right]{{\left[ x+1 \right]}_{n}}{{s}_{n}}\left( x \right)={{b}_{0}}{{\left( {{a}_{1}} \right)}^{n}}$, then in accordance with the Theorem 2.3.
$${{h}_{n}}\left( 1 \right)={{b}_{0}}{{\left( {{a}_{1}} \right)}^{n}}\frac{\left( 2n \right)!}{n!}.$$

If $b\left( x \right)=1$, $a\left( x \right)={{\left( 1-x \right)}^{-1}}$, then 
$${{h}_{n}}\left( x \right)=\left( n+1 \right)!{{N}_{n}}\left( x \right)={{\left( 1-x \right)}^{2n+1}}\sum\limits_{m=0}^{\infty }{{{\left[ m+1 \right]}_{n}}\left( \begin{matrix}
   m+n-1  \\
   n  \\
\end{matrix} \right)}{{x}^{m}},$$
where 
$${{N}_{0}}\left( x \right)=1,  \qquad{{N}_{n}}\left( x \right)=\frac{1}{n}\sum\limits_{m=0}^{n}{\left( \begin{matrix}
   n  \\
   m-1  \\
\end{matrix} \right)\left( \begin{matrix}
   n  \\
   n-m  \\
\end{matrix} \right){{x}^{m}}}$$ 
are the Narayana polynomials. In this connection, we will called polynomials ${{h}_{n}}\left( x \right)$ the generalized Narayana polynomials.

We introduce the matrices ${{F}_{n}}$:
$${{F}_{n}}=\frac{\left( 2n \right)!}{n!}{{U}_{2n}}\left( {{\left[ x+1 \right]}_{n}},x \right){{I}_{n}},  \qquad F_{n}^{-1}=\frac{n!}{\left( 2n \right)!}{{\left( {{\left[ x+1 \right]}_{n}},x \right)}^{-1}}U_{2n}^{-1}{{I}_{n}},$$
$${{F}_{n}}{{x}^{p}}=\frac{\left( 2n \right)!}{n!}{{U}_{2n}}{{x}^{p}}{{\left[ x+1 \right]}_{n}}={{\left( 1-x \right)}^{2n+1}}\sum\limits_{m=0}^{\infty }{{{m}^{p}}}\left( \begin{matrix}
   m+n  \\
   n  \\
\end{matrix} \right){{x}^{m}},$$
$$F_{n}^{-1}{{x}^{p}}=\frac{n!}{\left( 2n \right)!}{{\left( x \right)}_{p}}{{\left[ x+n+1 \right]}_{n-p}}.$$
For example,
$${{F}_{1}}=\left( \begin{matrix}
   1 & 0  \\
   -1 & 2  \\
\end{matrix} \right),  \qquad{{F}_{2}}=\left( \begin{matrix}
   1 & 0 & 0  \\
   -2 & 3 & 3  \\
   1 & -3 & 9  \\
\end{matrix} \right), \qquad{{F}_{3}}=\left( \begin{matrix}
   1 & 0 & 0 & 0  \\
   -3 & 4 & 4 & 4  \\
   3 & -8 & 12 & 52  \\
   -1 & 4 & -16 & 64  \\
\end{matrix} \right);$$
$$F_{1}^{-1}=\frac{1}{2!}\left( \begin{matrix}
   2 & 0  \\
   1 & 1  \\
\end{matrix} \right), \quad F_{2}^{-1}=\frac{2!}{4!}\left( \begin{matrix}
   12 & 0 & 0  \\
   7 & 3 & -1  \\
   1 & 1 & 1  \\
\end{matrix} \right), \quad F_{3}^{-1}=\frac{3!}{6!}\left( \begin{matrix}
   120 & 0 & 0 & 0  \\
   74 & 20 & -4 & 2  \\
   15 & 9 & 3 & -3  \\
   1 & 1 & 1 & 1  \\
\end{matrix} \right).$$
Then ${{F}_{n}}{{s}_{n}}\left( x \right)={{h}_{n}}\left( x \right)$.
\\{\bfseries Example 3.1.} We explain the identity
$${{F}_{n}}{{\left[ x+n+1 \right]}_{n}}=\frac{\left( 2n \right)!}{n!}.$$
Matrix ${{\left( {{\left( xa\left( x \right) \right)}^{\prime }},xa\left( x \right) \right)}_{E}}$ is similar to the matrix $\left( a\left( x \right),xa\left( x \right) \right)$ in the sense  that \linebreak $\left[ n,\searrow  \right]{{\left( {{\left( xa\left( x \right) \right)}^{\prime }},xa\left( x \right) \right)}_{E}}=\left( {1}/{x}\; \right)\left[ n,\searrow  \right]{{\left( 1,xa\left( x \right) \right)}_{E}}$, $n>0$. This is a consequence of the identities
$$\left[ {{x}^{n}} \right]\left( 1+x{{\left( \log a\left( x \right) \right)}^{\prime }} \right){{a}^{m}}\left( x \right)=\frac{m+n}{m}\left[ {{x}^{n}} \right]{{a}^{m}}\left( x \right),$$
$$\left[ {{x}^{n}} \right]{{\left( xa\left( x \right) \right)}^{\prime }}{{a}^{m}}\left( x \right)=\frac{m+n+1}{m+1}\left[ {{x}^{n}} \right]{{a}^{m+1}}\left( x \right).$$
It also shows that
$$\left[ n,\to  \right]{{\left( {{\left( xa\left( x \right) \right)}^{\prime }},\log a\left( x \right) \right)}_{E}}=\left( x+n+1 \right){{\tilde{u}}_{n}}\left( x+1 \right),  \qquad n>0,$$
where ${{\tilde{u}}_{n}}\left( x \right)=\left( {1}/{x}\; \right){{u}_{n}}\left( x \right)$. If $a\left( x \right)=C\left( x \right)$, where $C\left( x \right)$  is the Catalan series, then
$${{\tilde{u}}_{n}}\left( x \right)={{\left[ x+n+1 \right]}_{n-1}},  \qquad\left( x+n+1 \right){{\tilde{u}}_{n}}\left( x+1 \right)={{\left[ x+n+1 \right]}_{n}}.$$
Hence, numerator polynomials of the matrix ${{\left( 1,xC\left( x \right) \right)}_{E}}$ are the monomials $\left( {\left( 2n \right)!}/{n!}\; \right)x$. 
{\bfseries Theorem 3.1.}
$${{F}_{n}}{{E}^{n+1}}\left( 1,-x \right)F_{n}^{-1}={{\left( -1 \right)}^{n}}{{J}_{n}}.$$
{\bfseries Proof.}
$${{E}^{n+1}}\left( 1,-x \right){{\left( x \right)}_{p}}{{\left[ x+n+1 \right]}_{n-p}}={{\left( -x-n-1 \right)}_{p}}{{\left[ -x \right]}_{n-p}}={{\left( -1 \right)}^{n}}{{\left( x \right)}_{n-p}}{{\left[ x+n+1 \right]}_{p}},$$
or
$${{E}^{n+1}}\left( 1,-x \right)F_{n}^{-1}={{\left( -1 \right)}^{n}}F_{n}^{-1}{{J}_{n}}. \qquad\square $$

Denote ${{\left( 1,xa\left( x \right) \right)}^{-1}}=\left( 1,x\bar{a}\left( x \right) \right)$. 
\\{\bfseries Theorem 3.2.} \emph{ Polynomials ${{\left( -1 \right)}^{n}}{{J}_{n}}{{h}_{n}}\left( x \right)$ are the numerator polynomials of the matrix  ${{\left( b\left( x\bar{a}\left( x \right) \right){{\left( x\bar{a}\left( x \right) \right)}^{\prime }},x\bar{a}\left( x \right) \right)}_{E}}$.}
\\{\bfseries Proof.}  By the Theorem 3.1.
$${{\left( -1 \right)}^{n}}{{J}_{n}}{{h}_{n}}\left( x \right)={{F}_{n}}{{s}_{n}}\left( -x-n-1 \right).$$
We will imagine  an infinite table whose $k$th row, $k=0$, $\pm 1$, $\pm 2$, … , has the generating function
$$b\left( x \right){{a}^{k}}\left( x \right)=\sum\limits_{n=0}^{\infty }{\frac{{{s}_{n}}\left( k \right)}{n!}}{{x}^{n}}.$$
If the row numbers of the table decrease from top to bottom, then, according to the Lagrange inversion theorem, the $k$th descending diagonal of table has the generating function
$$b\left( x\bar{a}\left( x \right) \right)\left( 1+x{{\left( \log \bar{a}\left( x \right) \right)}^{\prime }} \right){{\left( \frac{1}{\bar{a}\left( x \right)} \right)}^{k}}=\sum\limits_{n=0}^{\infty }{\frac{{{s}_{n}}\left( k-n \right)}{n!}}{{x}^{n}}.$$
Hence,
$$b\left( x\bar{a}\left( x \right) \right){{\left( x\bar{a}\left( x \right) \right)}^{\prime }}{{\left( \bar{a}\left( x \right) \right)}^{k}}=\sum\limits_{n=0}^{\infty }{\frac{{{s}_{n}}\left( -k-n-1 \right)}{n!}}{{x}^{n}}. \qquad\square $$

Denote
$$\left[ n,\searrow  \right]{{\left( 1,xa\left( x \right) \right)}_{E}}=\frac{{{\varphi }_{n}}\left( x \right)}{{{\left( 1-x \right)}^{2n+1}}} , \qquad\left[ n,\searrow  \right]\left( 1,xa\left( x \right) \right)_{E}^{-1}=\frac{\varphi _{n}^{\left[ -1 \right]}\left( x \right)}{{{\left( 1-x \right)}^{2n+1}}}.$$
Then
$$\varphi _{n}^{\left[ -1 \right]}\left( x \right)={{\left( -1 \right)}^{n}}x{{J}_{n}}{{\varphi }_{n}}\left( x \right),\qquad n>0.$$
Hence, if the matrix $\left( 1,xa\left( x \right) \right)$ is the pseudo-involution [3], i.e., ${{\left( 1,xa\left( x \right) \right)}^{-1}}=\left( 1,xa\left( -x \right) \right)$,  then  ${{\varphi }_{n}}\left( x \right)=x{{J}_{n}}{{\varphi }_{n}}\left( x \right)$.

We will find a specific generating function for the sequence of polynomials ${{\varphi }_{n}}\left( t \right)$. Since
$$\frac{{{\varphi }_{n}}\left( t \right)}{\left( n+1 \right)!{{\left( 1-t \right)}^{2n+1}}}=\sum\limits_{m=0}^{\infty }{\frac{1}{n+1}}\left( \begin{matrix}
   n+m  \\
   m  \\
\end{matrix} \right)\left[ {{x}^{n}} \right]{{a}^{m}}\left( x \right){{t}^{m}}=$$
$$=\frac{1}{n+1}\left[ {{x}^{n}} \right]\frac{1}{{{\left( 1-ta\left( x \right) \right)}^{n+1}}}=\left[ {{x}^{n}} \right]b\left( x \right),$$
where according to the Lagrange inversion theorem
$$b\left( x \right)=\frac{1}{1-ta\left( xb\left( x \right) \right)} , \qquad\left( 1,xb\left( x \right) \right)={{\left( 1,x\left( 1-ta\left( x \right) \right) \right)}^{-1}},$$
then
$$\sum\limits_{n=0}^{\infty }{{{\varphi }_{n}}}\left( t \right)\frac{{{x}^{n}}}{\left( n+1 \right)!}=\left( 1-t \right)b\left( x{{\left( 1-t \right)}^{2}} \right).$$
{\bfseries Example 3.2.}
$$a\left( x \right)=\frac{1}{1-x},   \qquad b\left( x \right)=\frac{1+x-t-\sqrt{1-2x-2t-2xt+{{x}^{2}}+{{t}^{2}}}}{2x},$$
$$\sum\limits_{n=0}^{\infty }{{{\varphi }_{n}}\left( t \right)}\frac{{{x}^{n}}}{\left( n+1 \right)!}=\sum\limits_{n=0}^{\infty }{{{N}_{n}}}\left( t \right){{x}^{n}}=\frac{1+x\left( 1-t \right)-\sqrt{1-2x\left( 1+t \right)+{{x}^{2}}{{\left( 1-t \right)}^{2}}}}{2x}.$$

Generating functions of the second-order Euler numbers which are the coefficients of numerator polynomials of the matrices ${{\left( 1,{{e}^{x}}-1 \right)}_{E}}$, ${{\left( {{e}^{x}},{{e}^{x}}-1 \right)}_{E}}$ are considered in [17]. 
\section{Сonnection matrices }
We introduce the matrices ${{S}_{n}}={{F}_{n}}U_{n}^{-1}$,  $S_{n}^{-1}={{U}_{n}}F_{n}^{-1}$ . For example,
$${{S}_{1}}=\left( \begin{matrix}
   1 & 0  \\
   1 & 2  \\
\end{matrix} \right),  \qquad{{S}_{2}}=2!\left( \begin{matrix}
   1 & 0 & 0  \\
   4 & 3 & 0  \\
   1 & 3 & 6  \\
\end{matrix} \right), \qquad{{S}_{3}}=3!\left( \begin{matrix}
   1 & 0 & 0 & 0  \\
   9 & 4 & 0 & 0  \\
   9 & 12 & 10 & 0  \\
   1 & 4 & 10 & 20  \\
\end{matrix} \right),$$
$${{S}_{4}}=4!\left( \begin{matrix}
   1 & 0 & 0 & 0 & 0  \\
   16 & 5 & 0 & 0 & 0  \\
   36 & 30 & 15 & 0 & 0  \\
   16 & 30 & 40 & 35 & 0  \\
   1 & 5 & 15 & 35 & 70  \\
\end{matrix} \right);$$
$$S_{1}^{-1}=\frac{1}{2!}\left( \begin{matrix}
   2 & 0  \\
   -1 & 1  \\
\end{matrix} \right), \qquad S_{2}^{-1}=\frac{2!}{4!}\left( \begin{matrix}
   6 & 0 & 0  \\
   -8 & 2 & 0  \\
   3 & -1 & 1  \\
\end{matrix} \right), \qquad S_{3}^{-1}=\frac{3!}{6!}\left( \begin{matrix}
   20 & 0 & 0 & 0  \\
   -45 & 5 & 0 & 0  \\
   36 & -6 & 2 & 0  \\
   -10 & 2 & -1 & 1  \\
\end{matrix} \right),$$
$$S_{4}^{-1}=\frac{4!}{8!}\left( \begin{matrix}
   70 & 0 & 0 & 0 & 0  \\
   -224 & 14 & 0 & 0 & 0  \\
   280 & -28 & {14}/{3}\; & 0 & 0  \\
   -160 & 20 & {-16}/{3}\; & 2 & 0  \\
   35 & -5 & {5}/{3}\; & -1 & 1  \\
\end{matrix} \right).$$
Then ${{S}_{n}}{{g}_{n}}\left( x \right)={{h}_{n}}\left( x \right)$.
\\{\bfseries Theorem 4.1.} 
$${{S}_{n}}=V_{n}^{-1}{{C}_{n}}{{V}_{n}},  \qquad{{C}_{n}}{{x}^{p}}=\frac{\left( n+p \right)!}{p!}{{x}^{p}}.$$
{\bfseries Proof.} We use  the Theorem 2.4. and the identities
$$U_{n}^{-1}V_{n}^{-1}{{x}^{p}}=\frac{n!}{p!}{{\left( x \right)}_{p}},  \qquad U_{n+p}^{-1}{{x}^{p}}={{\left( x \right)}_{p}}{{\left[ x+1 \right]}_{n}}.$$
Then
$${{F}_{n}}U_{n}^{-1}V_{n}^{-1}{{x}^{p}}=\frac{\left( 2n \right)!}{n!}{{U}_{2n}}\frac{n!}{p!}{{\left( x \right)}_{p}}{{\left[ x+1 \right]}_{n}}=$$
$$=\frac{\left( 2n \right)!}{p!}{{\left( 1-x \right)}^{n-p}}\frac{\left( n+p \right)!}{\left( 2n \right)!}{{U}_{n+p}}{{\left( x \right)}_{p}}{{\left[ x+1 \right]}_{n}}=\frac{\left( n+p \right)!}{p!}{{\left( 1-x \right)}^{n-p}}{{x}^{p}},$$
or
$${{S}_{n}}V_{n}^{-1}=V_{n}^{-1}{{C}_{n}}. \qquad\square $$

For example,
$${{S}_{3}}=3!\left( \begin{matrix}
   1 & 0 & 0 & 0  \\
   -3 & 1 & 0 & 0  \\
   3 & -2 & 1 & 0  \\
   -1 & 1 & -1 & 1  \\
\end{matrix} \right)\left( \begin{matrix}
   1 & 0 & 0 & 0  \\
   0 & 4 & 0 & 0  \\
   0 & 0 & 10 & 0  \\
   0 & 0 & 0 & 20  \\
\end{matrix} \right)\left( \begin{matrix}
   1 & 0 & 0 & 0  \\
   3 & 1 & 0 & 0  \\
   3 & 2 & 1 & 0  \\
   1 & 1 & 1 & 1  \\
\end{matrix} \right).$$
{\bfseries Theorem 4.2.} 
$${{S}_{n}}{{x}^{p}}=\frac{\left( n+p \right)!\left( n-p \right)!}{n!}\sum\limits_{m=p}^{n}{\left( \begin{matrix}
   n  \\
   m-p  \\
\end{matrix} \right)\left( \begin{matrix}
   n  \\
   n-m  \\
\end{matrix} \right){{x}^{m}}}.$$
{\bfseries Proof.} 
$$\left[ m,\to  \right]V_{n}^{-1}=\sum\limits_{i=0}^{m}{\left( \begin{matrix}
   m-n-1  \\
   m-i  \\
\end{matrix} \right){{x}^{i}}}=\sum\limits_{i=0}^{m}{{{\left( -1 \right)}^{m-i}}}\left( \begin{matrix}
   n-i  \\
   m-i  \\
\end{matrix} \right){{x}^{i}},$$
$${{C}_{n}}{{V}_{n}}{{x}^{p}}=\sum\limits_{i=p}^{n}{\frac{\left( n+i \right)!}{i!}\left( \begin{matrix}
   n-p  \\
   i-p  \\
\end{matrix} \right)}{{x}^{i}},$$
$$\left[ {{x}^{m}} \right]V_{n}^{-1}{{C}_{n}}{{V}_{n}}{{x}^{p}}=\sum\limits_{i=p}^{m}{{{\left( -1 \right)}^{m-i}}\left( \begin{matrix}
   n-i  \\
   m-i  \\
\end{matrix} \right)\frac{\left( n+i \right)!}{i!}\left( \begin{matrix}
   n-p  \\
   i-p  \\
\end{matrix} \right)}=$$
$$=\frac{\left( n+p \right)!\left( n-p \right)!}{\left( n-m \right)!m!}\sum\limits_{i=p}^{m}{{{\left( -1 \right)}^{m-i}}}\left( \begin{matrix}
   n+i  \\
   i-p  \\
\end{matrix} \right)\left( \begin{matrix}
   m  \\
   i  \\
\end{matrix} \right)=$$ 
$$=\frac{\left( n+p \right)!\left( n-p \right)!}{\left( n-m \right)!m!}{{\left( -1 \right)}^{m-p}}\left( \begin{matrix}
   m-n-p-1  \\
   m-p  \\
\end{matrix} \right)=
\frac{\left( n+p \right)!\left( n-p \right)!}{n!}\left( \begin{matrix}
   n  \\
   m-p  \\
\end{matrix} \right)\left( \begin{matrix}
   n  \\
   n-m  \\
\end{matrix} \right). \,\square $$
{\bfseries Theorem 4.3.} 
$$S_{n}^{-1}{{x}^{p}}=\frac{p!\left( n-p \right)!}{\left( 2n \right)!}\sum\limits_{m=p}^{n}{\left( \begin{matrix}
   -n  \\
   m-p  \\
\end{matrix} \right)\left( \begin{matrix}
   2n  \\
   n-m  \\
\end{matrix} \right){{x}^{m}}}.$$
{\bfseries Proof.} 
$$\left[ {{x}^{m}} \right]V_{n}^{-1}C_{n}^{-1}{{V}_{n}}{{x}^{p}}=\sum\limits_{i=p}^{m}{{{\left( -1 \right)}^{m-i}}\left( \begin{matrix}
   n-i  \\
   m-i  \\
\end{matrix} \right)\frac{i!}{\left( n+i \right)!}\left( \begin{matrix}
   n-p  \\
   i-p  \\
\end{matrix} \right)}=$$
$$=\frac{p!\left( n-p \right)!}{\left( n-m \right)!\left( n+m \right)!}\sum\limits_{i=p}^{m}{{{\left( -1 \right)}^{m-i}}}\left( \begin{matrix}
   i  \\
   i-p  \\
\end{matrix} \right)\left( \begin{matrix}
   n+m  \\
   m-i  \\
\end{matrix} \right)=$$ 
$$=\frac{p!\left( n-p \right)!}{\left( n-m \right)!\left( n+m \right)!}{{\left( -1 \right)}^{m-p}}\left( \begin{matrix}
   n+m-p-1  \\
   m-p  \\
\end{matrix} \right)=
\frac{p!\left( n-p \right)!}{\left( 2n \right)!}\left( \begin{matrix}
   -n  \\
   m-p  \\
\end{matrix} \right)\left( \begin{matrix}
   2n  \\
   n-m  \\
\end{matrix} \right). \,\square $$
{\bfseries Example 4.1.} Numerator polynomials of the Pascal matrix are equal to one.  Hence, numerator polynomials of the matrix ${{\left| {{e}^{x}} \right|}^{-1}}P\left| {{e}^{x}} \right|$ are the polynomials ${{S}_{n}}{{x}^{0}}=n!{}^{B}{{N}_{n}}\left( x \right)$, where ${}^{B}{{N}_{n}}\left( x \right)$ are the so-called Narayana polynomials of type $B$:
$${}^{B}{{N}_{n}}\left( x \right)=\sum\limits_{m=0}^{n}{{{\left( \begin{matrix}
   n  \\
   m  \\
\end{matrix} \right)}^{2}}{{x}^{m}}}={{\left( 1-x \right)}^{2n+1}}\sum\limits_{m=0}^{\infty }{{{\left( \begin{matrix}
   m+n  \\
   n  \\
\end{matrix} \right)}^{2}}{{x}^{m}}}.$$

Pascal matrix belongs to a subgroup of the Riordan group formed by the matrices of the form $\left( a\left( x \right),xa\left( x \right) \right)$ and called the Bell subgroup. Numerator polynomials with the property ${{h}_{n}}\left( x \right)={{J}_{n}}{{h}_{n}}\left( x \right)$ will be called symmetric.
\\{\bfseries Theorem 4.4.} \emph{In the exponential Bell subgroup, only the matrices ${{\left| {{e}^{x}} \right|}^{-1}}{{P}^{\varphi }}\left| {{e}^{x}} \right|$ have the symmetric numerator polynomials.}
\\{\bfseries Proof.} Let ${{h}_{n}}\left( x \right)$ are the numerator polynomials of the matrix ${{\left( a\left( x \right),xa\left( x \right) \right)}_{E}}$. Since $a\left( x\bar{a}\left( x \right) \right)={1}/{\bar{a}\left( x \right)}\;$, then by the Theorem 3.2. the polynomials ${{\left( -1 \right)}^{n}}{{J}_{n}}{{h}_{n}}\left( x \right)$ are the numerator polynomials of the matrix
$${{\left( 1+x{{\left( \log \bar{a}\left( x \right) \right)}^{\prime }},x\bar{a}\left( x \right) \right)}_{E}}=\left( 1+x{{\left( \log a\left( x \right) \right)}^{\prime }},xa\left( x \right) \right)_{E}^{-1}.$$
Thus, the necessary condition for the symmetry of numerator polynomials is the condition $a\left( x \right)=1+x{{\left( \log a\left( x \right) \right)}^{\prime }}$. It comes down to the condition ${{a}^{2}}\left( x \right)={{\left( xa\left( x \right) \right)}^{\prime }}$, or $\sum\nolimits_{m=0}^{n}{{{a}_{n-m}}{{a}_{m}}}=\left( n+1 \right){{a}_{n}}$, that is feasible only in the case ${{a}_{n}}={{\left( {{a}_{1}} \right)}^{n}}. \qquad\square $

$n$th numerator polynomial of the matrix $A$ (i.e., ordinary or exponential Riordan matrix) will be denoted by $\left[ {{P}_{n}}\left( x \right) \right]A$. 
\\{\bfseries Example 4.2.} Since $\left[ {{P}_{n}}\left( x \right) \right]\left( 1+x,x\left( 1+x \right) \right)={{x}^{n-1}}$, $n>0$, then
$$\left[ {{P}_{n}}\left( x \right) \right]{{\left( 1+x,x\left( 1+x \right) \right)}_{E}}={{S}_{n}}{{x}^{n-1}}=\frac{\left( 2n \right)!}{n!2}\left( 1+x \right){{x}^{n-1}}:$$
$${{\left( 1+x,x\left( 1+x \right) \right)}_{E}}={{\left| {{e}^{x}} \right|}^{-1}}\left( \begin{matrix}
   1 & 0 & 0 & 0 & 0 & 0 & \cdots   \\
   1 & 1 & 0 & 0 & 0 & 0 & \cdots   \\
   0 & 2 & 1 & 0 & 0 & 0 & \cdots   \\
   0 & 1 & 3 & 1 & 0 & 0 & \cdots   \\
   0 & 0 & 3 & 4 & 1 & 0 & \cdots   \\
   0 & 0 & 1 & 6 & 5 & 1 & \cdots   \\
   \vdots  & \vdots  & \vdots  & \vdots  & \vdots  & \vdots  & \ddots   \\
\end{matrix} \right)\left| {{e}^{x}} \right|.$$
Since
$$\left( 1+x{{\left( \log \left( 1+x \right) \right)}^{\prime }},x\left( 1+x \right) \right)_{E}^{-1}={{\left( 1+x{{\left( \log C\left( -x \right) \right)}^{\prime }},xC\left( -x \right) \right)}_{E}},$$
then
$$\left[ {{P}_{n}}\left( x \right) \right]{{\left( 1+x{{\left( \log C\left( x \right) \right)}^{\prime }},xC\left( x \right) \right)}_{E}}=\frac{\left( 2n \right)!}{n!2}\left( 1+x \right):$$
$${{\left( 1+x{{\left( \log C\left( x \right) \right)}^{\prime }},xC\left( x \right) \right)}_{E}}={{\left| {{e}^{x}} \right|}^{-1}}\left( \begin{matrix}
   1 & 0 & 0 & 0 & 0 & 0 & \cdots   \\
   1 & 1 & 0 & 0 & 0 & 0 & \cdots   \\
   3 & 2 & 1 & 0 & 0 & 0 & \cdots   \\
   10 & 6 & 3 & 1 & 0 & 0 & \cdots   \\
   35 & 20 & 10 & 4 & 1 & 0 & \cdots   \\
   126 & 70 & 35 & 15 & 5 & 1 & \cdots   \\
   \vdots  & \vdots  & \vdots  & \vdots  & \vdots  & \vdots  & \ddots   \\
\end{matrix} \right)\left| {{e}^{x}} \right|.$$
{\bfseries Example 4.3.} Since $\left[ {{P}_{n}}\left( x \right) \right]{{\left( {{\left( xC\left( x \right) \right)}^{\prime }},xC\left( x \right) \right)}_{E}}={\left( 2n \right)!}/{n!}\;$, then
$$\left[ {{P}_{n}}\left( x \right) \right]\left( {{\left( xC\left( x \right) \right)}^{\prime }},xC\left( x \right) \right)=\frac{\left( 2n \right)!}{n!}S_{n}^{-1}{{x}^{0}}=\sum\limits_{m=0}^{n}{\left( \begin{matrix}
   -n  \\
   m  \\
\end{matrix} \right)\left( \begin{matrix}
   2n  \\
   n-m  \\
\end{matrix} \right){{x}^{m}}}:$$
$$\left( {{\left( xC\left( x \right) \right)}^{\prime }},xC\left( x \right) \right)=\left( \begin{matrix}
   1 & 0 & 0 & 0 & 0 & 0 & \cdots   \\
   2 & 1 & 0 & 0 & 0 & 0 & \cdots   \\
   6 & 3 & 1 & 0 & 0 & 0 & \cdots   \\
   20 & 10 & 4 & 1 & 0 & 0 & \cdots   \\
   70 & 35 & 15 & 5 & 1 & 0 & \cdots   \\
   252 & 126 & 56 & 21 & 6 & 1 & \cdots   \\
   \vdots  & \vdots  & \vdots  & \vdots  & \vdots  & \vdots  & \ddots   \\
\end{matrix} \right).$$
Respectively,
$${{\left( -1 \right)}^{n}}\sum\limits_{m=0}^{n}{\left( \begin{matrix}
   2n  \\
   m  \\
\end{matrix} \right)\left( \begin{matrix}
   -n  \\
   n-m  \\
\end{matrix} \right){{x}^{m}}}=\left[ {{P}_{n}}\left( x \right) \right]\left( 1+x{{\left( \operatorname{logC}\left( x \right) \right)}^{\prime }},x{{C}^{-1}}\left( x \right) \right):$$
$$\left( 1+x{{\left( \operatorname{logC}\left( x \right) \right)}^{\prime }},x{{C}^{-1}}\left( x \right) \right)=\left( \begin{matrix}
   1 & 0 & 0 & 0 & 0 & 0 & \cdots   \\
   1 & 1 & 0 & 0 & 0 & 0 & \cdots   \\
   3 & 0 & 1 & 0 & 0 & 0 & \cdots   \\
   10 & 1 & -1 & 1 & 0 & 0 & \cdots   \\
   35 & 4 & 0 & -2 & 1 & 0 & \cdots   \\
   126 & 15 & 1 & 0 & -3 & 1 & \cdots   \\
   \vdots  & \vdots  & \vdots  & \vdots  & \vdots  & \vdots  & \ddots   \\
\end{matrix} \right).$$

Let $_{\left( \beta  \right)}a\left( x \right)$ is the generalized binomial series. Denote
$$[n,\searrow ]{{\left( 1,{{x}_{\left( \beta  \right)}}a\left( x \right) \right)}_{E}}=\frac{_{\left( \beta  \right)}{{\varphi }_{n}}\left( x \right)}{{{\left( 1-x \right)}^{2n+1}}}=\sum\limits_{m=0}^{\infty }{\frac{m}{m+\beta n}}\left( \begin{matrix}
   m+\beta n  \\
   n  \\
\end{matrix} \right){{\left[ m+1 \right]}_{n}}{{x}^{m}}.$$
{\bfseries Theorem 4.5.} 
$$_{\left( \beta  \right)}{{\varphi }_{n}}\left( x \right)=\frac{\left( n+1 \right)!}{n}\sum\limits_{m=0}^{n}{\left( \begin{matrix}
   n\left( 2-\beta  \right)  \\
   m-1  \\
\end{matrix} \right)\left( \begin{matrix}
   n\beta   \\
   n-m  \\
\end{matrix} \right){{x}^{m}}}, \qquad n>0.\eqno (3)$$
{\bfseries Proof.} We prove the theorem by analogy with the proof of  Theorem 2.5.
$$\left[ {{x}^{m}} \right]{}_{\left( \beta  \right)}{{\varphi }_{n}}\left( x \right)=\left[ {{x}^{m}} \right]{{S}_{n}}{}_{\left( \beta  \right)}{{\alpha }_{n}}\left( x \right)=\left[ {{x}^{m}} \right]V_{n}^{-1}{{C}_{n}}{}_{\left( \beta  \right)}{{v}_{n}}\left( x \right)=$$
$$=\sum\limits_{i=0}^{m}{{{\left( -1 \right)}^{m-i}}}\left( \begin{matrix}
   n-i  \\
   m-i  \\
\end{matrix} \right)\frac{\left( n+i \right)!}{i!}\frac{i}{n}\left( \begin{matrix}
   n\beta   \\
   n-i  \\
\end{matrix} \right)\frac{\left( n\beta -n+m \right)!}{\left( n\beta -n+m \right)!}=$$
$$=\frac{\left( n+1 \right)!}{n}\left( \begin{matrix}
   n\beta   \\
   n-m  \\
\end{matrix} \right)\sum\limits_{i=0}^{m}{{{\left( -1 \right)}^{m-i}}\left( \begin{matrix}
   n+i  \\
   i-1  \\
\end{matrix} \right)}\left( \begin{matrix}
   n\beta -n+m  \\
   m-i  \\
\end{matrix} \right)=$$ 
$$=\frac{\left( n+1 \right)!}{n}\left( \begin{matrix}
   n\beta   \\
   n-m  \\
\end{matrix} \right){{\left( -1 \right)}^{m-1}}\left( \begin{matrix}
   n\beta -2n+m-2  \\
   m-1  \\
\end{matrix} \right)=
\frac{\left( n+1 \right)!}{n}\left( \begin{matrix}
   n\beta   \\
   n-m  \\
\end{matrix} \right)\left( \begin{matrix}
   n\left( 2-\beta  \right)  \\
   m-1  \\
\end{matrix} \right). \,\square $$

Note that
$$_{\left( 0 \right)}{{\varphi }_{n}}\left( x \right)=\frac{\left( 2n \right)!}{n!}{{x}^{n}},  \qquad_{\left( 1 \right)}{{\varphi }_{n}}\left( x \right)=\left( n+1 \right)!{{N}_{n}}\left( x \right),   \qquad_{\left( 2 \right)}{{\varphi }_{n}}\left( x \right)=\frac{\left( 2n \right)!}{n!}x.$$
Since ${{\left( 1,x{}_{\left( \beta  \right)}a\left( x \right) \right)}^{-1}}=\left( 1,{{x}_{\left( \beta -1 \right)}}{{a}^{-1}}\left( x \right) \right)$,  ${}_{\left( \beta -1 \right)}{{a}^{-1}}\left( -x \right){{=}_{\left( 2-\beta  \right)}}a\left( x \right)$, then $_{\left( 2-\beta  \right)}{{\varphi }_{n}}\left( x \right)=x{{J}_{n}}{}_{\left( \beta  \right)}{{\varphi }_{n}}\left( x \right)$.
\section{ Generalized Lagrange series}
Let $_{\left( \beta  \right)}a\left( x \right)$ is the generalized binomial series. We agree to denote $_{\left( 0 \right)}a\left( x \right)=a\left( x \right)$. Then ${{\left( 1,x{{a}^{-\beta }}\left( x \right) \right)}^{-1}}=\left( 1,x{}_{\left( \beta  \right)}{{a}^{\beta }}\left( x \right) \right)$, ${{\left( 1,x{}_{\left( \beta  \right)}{{a}^{\varphi }}\left( x \right) \right)}^{-1}}=\left( 1,{{x}_{\left( \beta -\varphi  \right)}}{{a}^{-\varphi }}\left( x \right) \right)$. We extend the generalization that underlies this construction to any formal power series $a\left( x \right)$, ${{a}_{0}}=1$.  Let ${{\left( 1,x{{a}^{-1}}\left( x \right) \right)}^{-1}}=\left( 1,xb\left( x \right) \right)$. Denote $\left[ {{x}^{n}} \right]{{a}^{m}}\left( x \right)=a_{n}^{\left( m \right)}$. For Riordan matrices, two forms of the Lagrange inversion formula
$$f\left( x \right)={{f}_{0}}+\sum\limits_{n=1}^{\infty }{\frac{{{x}^{n}}}{{{a}^{n}}\left( x \right)}}\frac{1}{n}\left[ {{x}^{n-1}} \right]{f}'\left( x \right){{a}^{n}}\left( x \right),$$
$$\frac{f\left( x \right)}{1-x{{\left( \log a\left( x \right) \right)}^{\prime }}}=\sum\limits_{n=0}^{\infty }{\frac{{{x}^{n}}}{{{a}^{n}}\left( x \right)}}\left[ {{x}^{n}} \right]f\left( x \right){{a}^{n}}\left( x \right),$$
 where $f\left( x \right)$ is an arbitrary series, take the form
$$\left[ n,\to  \right]{{\left( 1,x{{a}^{-1}}\left( x \right) \right)}^{-1}}=\sum\limits_{m=0}^{n}{\frac{m}{n}a_{n-m}^{\left( n \right)}}{{x}^{m}}, \qquad n>0;$$
$$\left[ n,\to  \right]{{\left( 1-x{{\left( \log a\left( x \right) \right)}^{\prime }},x{{a}^{-1}}\left( x \right) \right)}^{-1}}=\sum\limits_{m=0}^{n}{a_{n-m}^{\left( n \right)}}{{x}^{m}}.$$
Hence,
 $$\left[ {{x}^{n}} \right]{{b}^{m}}\left( x \right)=\frac{m}{m+n}\left[ {{x}^{n}} \right]{{a}^{m+n}}\left( x \right)=\left[ {{x}^{n}} \right]\left( 1-x{{\left( \log a\left( x \right) \right)}^{\prime }} \right){{a}^{m+n}}\left( x \right),$$
$$\left[ {{x}^{n}} \right]\left( 1+x{{\left( \log b\left( x \right) \right)}^{\prime }} \right){{b}^{m}}\left( x \right)=\frac{m+n}{m}\left[ {{x}^{n}} \right]{{b}^{m}}\left( x \right)=\left[ {{x}^{n}} \right]{{a}^{m+n}}\left( x \right),$$
$${{\left( 1-x{{\left( \log a\left( x \right) \right)}^{\prime }},x{{a}^{-1}}\left( x \right) \right)}^{-1}}=\left( 1+x{{\left( \log b\left( x \right) \right)}^{\prime }},xb\left( x \right) \right).$$
Denote ${{\left( 1,x{{a}^{-\beta }}\left( x \right) \right)}^{-1}}=\left( 1,x{}_{\left( \beta  \right)}{{a}^{\beta }}\left( x \right) \right)$. Then
$$\left[ {{x}^{n}} \right]{}_{\left( \beta  \right)}{{a}^{\beta m}}\left( x \right)=\frac{\beta m}{\beta m+\beta n}\left[ {{x}^{n}} \right]{{a}^{\beta m+\beta n}}\left( x \right),$$ 
$$_{\left( \beta  \right)}{{a}^{\varphi }}\left( x \right)=\sum\limits_{n=0}^{\infty }{\frac{\varphi }{\varphi +\beta n}}\frac{{{u}_{n}}\left( \varphi +\beta n \right)}{n!}{{x}^{n}},  \qquad{{u}_{n}}\left( x \right)=\left[ n,\to  \right]{{\left( 1,\log a\left( x \right) \right)}_{E}}.$$
  Let ${{\left( 1,\log a\left( x \right) \right)}^{-1}}=\left( 1,q\left( x \right) \right)$. Then 

$$\left( 1,\log {}_{\left( \beta  \right)}a\left( x \right) \right)=\left( 1,{{x}_{\left( \beta  \right)}}{{a}^{\beta }}\left( x \right) \right)\left( 1,\log a\left( x \right) \right),$$
$${{\left( 1,\log {}_{\left( \beta  \right)}a\left( x \right) \right)}^{-1}}=\left( 1,q\left( x \right) \right)\left( 1,x{{a}^{-\beta }}\left( x \right) \right)=\left( 1,q\left( x \right){{e}^{-\beta x}} \right).$$

Denote.
$${}_{\left( \beta  \right)}{{u}_{n}}\left( x \right)=\left[ n,\to  \right]{{\left( 1,\log {}_{\left( \beta  \right)}a\left( x \right) \right)}_{E}},  \qquad_{\left( \beta  \right)}{{q}_{n}}\left( x \right)=\left( 1,\log {}_{\left( \beta  \right)}a\left( x \right) \right)_{E}^{-1}{{x}^{n}}.$$
Since ${{\left( 1,q\left( x \right){{e}^{-\beta x}} \right)}_{E}}{{x}^{n}}={{\left( {{e}^{-\beta nx}},q\left( x \right) \right)}_{E}}{{x}^{n}}$, then the system of identities take place:
$$_{\left( \beta  \right)}{{u}_{n}}\left( x \right)=\frac{x}{x+n\beta }{{u}_{n}}\left( x+n\beta  \right), \qquad_{\left( \beta  \right)}{{q}_{n}}\left( x \right)=\frac{1}{1+n\beta x}{{q}_{n}}\left( \frac{x}{1+n\beta x} \right),$$
$$\sum\limits_{n=0}^{\infty }{_{\left( \beta  \right)}{{u}_{n}}}\left( \varphi  \right){}_{\left( \beta  \right)}{{q}_{n}}\left( x \right)=\frac{1}{1-\varphi x}.$$

Series $_{\left( \beta  \right)}a\left( x \right)$ for integers $\beta =k$ (denoted by ${{S}_{k}}\left( x \right)$) were introduced in [18]. In [19], these series, called generalized Lagrange series, are considered in connection with the Riordan matrices. They have the following visual interpretation. We will imagine  the table ${{\left\{ b\left( x \right),{{a}^{\varphi }}\left( x \right) \right\}}_{0}}$ whose $k$th row, $k=0$, $\pm 1$, $\pm 2$, … , has the generating function $b\left( x \right){{a}^{\varphi k}}\left( x \right)$. Lagrange inversion formula connects the rows of the table with its ascending and descending diagonals in a certain way. For definiteness, we assume that $\varphi >0$ and the table row numbers decrease from top to bottom. We replace the $k$th row of the table ${{\left\{ b\left( x \right),{{a}^{\varphi }}\left( x \right) \right\}}_{0}}$ with the $k$th ascending diagonal, the resulting table will be denoted by ${{\left\{ b\left( x \right),{{a}^{\varphi }}\left( x \right) \right\}}_{1}}$. We perform the same operation with the table ${{\left\{ b\left( x \right),{{a}^{\varphi }}\left( x \right) \right\}}_{1}}$, the result will be denoted by  ${{\left\{ b\left( x \right),{{a}^{\varphi }}\left( x \right) \right\}}_{2}}$, etc. Tables obtained in a similar way for descending diagonals are numbered with negative numbers. Then it turns out that the $k$th row of the table ${{\left\{ b\left( x \right),{{a}^{\varphi }}\left( x \right) \right\}}_{v}}$ has the generating function 
$$b\left( {{x}_{\left( v\varphi  \right)}}{{a}^{v\varphi }}\left( x \right) \right)\left( 1+x{{\left( \log {}_{\left( v\varphi  \right)}{{a}^{v\varphi }}\left( x \right) \right)}^{\prime }} \right){}_{\left( v\varphi  \right)}{{a}^{\varphi k}}\left( x \right),$$
that can be written as the pair of mutually inverse identities:
$${{\left\{ b\left( x \right),{{a}^{\varphi }}\left( x \right) \right\}}_{v}}={{\left\{ b\left( {{x}_{\left( v\varphi  \right)}}{{a}^{v\varphi }}\left( x \right) \right)\left( 1+x{{\left( \log {}_{\left( v\varphi  \right)}{{a}^{v\varphi }}\left( x \right) \right)}^{\prime }} \right),{}_{\left( v\varphi  \right)}{{a}^{\varphi }}\left( x \right) \right\}}_{0}},$$
$${{\left\{ b\left( {{x}_{\left( v\varphi  \right)}}{{a}^{v\varphi }}\left( x \right) \right)\left( 1+x{{\left( \log {}_{\left( v\varphi  \right)}{{a}^{v\varphi }}\left( x \right) \right)}^{\prime }} \right),{}_{\left( v\varphi  \right)}{{a}^{\varphi }}\left( x \right) \right\}}_{-v}}={{\left\{ b\left( x \right),{{a}^{\varphi }}\left( x \right) \right\}}_{0}}.$$
Operators mapping rows of one table to rows of another, and vice versa, correspond to the matrices
$$\left( 1+x{{\left( \log {}_{\left( \beta  \right)}{{a}^{\beta }}\left( x \right) \right)}^{\prime }},x{}_{\left( \beta  \right)}{{a}^{\beta }}\left( x \right) \right), \quad\left( 1+x{{\left( \log {{a}^{-\beta }}\left( x \right) \right)}^{\prime }},x{{a}^{-\beta }}\left( x \right) \right), \quad\beta =v\varphi.$$

In the Riordan matrices theory, generalized Lagrange series do not have a generally accepted notation and are present implicitly in the form of various constructions. These constructions are presented in the papers [20] – [26].
\section{Polynomials $_{\left( \beta  \right)}{{g}_{n}}\left( x \right)$. }
Let $\left[ n,\to  \right]{{\left( b\left( x \right),\log a\left( x \right) \right)}_{E}}={{s}_{n}}\left( x \right)$. Then by definition of the series $_{\left( \beta  \right)}a\left( x \right)$
$$\left[ n,\to  \right]{{\left( b\left( x{}_{\left( \beta  \right)}{{a}^{\beta }}\left( x \right) \right)\left( 1+x{{\left( \log {}_{\left( \beta  \right)}{{a}^{\beta }}\left( x \right) \right)}^{\prime }} \right),\log {}_{\left( \beta  \right)}a\left( x \right) \right)}_{E}}={{s}_{n}}\left( x+\beta n \right).$$
Denote 
$$\left[ {{P}_{n}}\left( x \right) \right]\left( b\left( x{}_{\left( \beta  \right)}{{a}^{\beta }}\left( x \right) \right)\left( 1+x{{\left( \log {}_{\left( \beta  \right)}{{a}^{\beta }}\left( x \right) \right)}^{\prime }} \right),x{}_{\left( \beta  \right)}a\left( x \right) \right){{=}_{\left( \beta  \right)}}{{g}_{n}}\left( x \right).$$
We introduce the matrices $G_{n}^{\beta }={{U}_{n}}{{E}^{n\beta }}U_{n}^{-1}$. Then $G_{n}^{\beta }{{g}_{n}}\left( x \right)={}_{\left( \beta  \right)}{{g}_{n}}\left( x \right)$. For example,
$${{G}_{1}}=\left( \begin{matrix}
   2 & 1  \\
   -1 & 0  \\
\end{matrix} \right), \qquad{{G}_{2}}=\left( \begin{matrix}
   6 & 3 & 1  \\
   -8 & -3 & 0  \\
   3 & 1 & 0  \\
\end{matrix} \right),   \qquad{{G}_{3}}=\left( \begin{matrix}
   20 & 10 & 4 & 1  \\
   -45 & -20 & -6 & 0  \\
   36 & 15 & 4 & 0  \\
   -10 & -4 & -1 & 0  \\
\end{matrix} \right),$$
Since $\left[ x \right]{}_{\left( \beta  \right)}a\left( x \right)={{a}_{1}}$, then $_{\left( \beta  \right)}{{g}_{n}}\left( 1 \right)={{g}_{n}}\left( 1 \right)$; therefore, the sum of elements of each column of the matrix $G_{n}^{\beta }$ is equal to one.
\\{\bfseries Theorem 6.1.} \emph{}
$$G_{n}^{-\beta }={{J}_{n}}G_{n}^{\beta }{{J}_{n}}.$$
{\bfseries Proof.}  Since ${{E}^{\varphi }}\left( 1,-x \right)=\left( 1,-x \right){{E}^{-\varphi }}$, then by the Theorem 2.1.
$${{J}_{n}}{{U}_{n}}{{E}^{n\beta }}U_{n}^{-1}{{J}_{n}}={{U}_{n}}E\left( 1,-x \right){{E}^{n\beta }}E\left( 1,-x \right)U_{n}^{-1}={{U}_{n}}{{E}^{-n\beta }}U_{n}^{-1}.\qquad\square $$ 

Thus,
$$G_{1}^{-1}=\left( \begin{matrix}
   0 & -1  \\
   1 & 2  \\
\end{matrix} \right), \qquad G_{2}^{-1}=\left( \begin{matrix}
   0 & 1 & 3  \\
   0 & -3 & -8  \\
   1 & 3 & 6  \\
\end{matrix} \right), \qquad G_{3}^{-1}=\left( \begin{matrix}
   0 & -1 & -4 & -10  \\
   0 & 4 & 15 & 36  \\
   0 & -6 & -20 & -45  \\
   1 & 4 & 10 & 20  \\
\end{matrix} \right).$$
{\bfseries Theorem 6.2.} \emph{}
$$G_{n}^{\beta }=V_{n}^{-1}{{\left( {{\left( 1+x \right)}^{n\beta }},x \right)}^{T}}{{V}_{n}}.$$
{\bfseries Proof.} Since
$$n!\left| {{e}^{x}} \right|{{V}_{n}}{{U}_{n}}{{x}^{p}}=\left[ p,\to  \right]{{\left( 1,{{e}^{x}}-1 \right)}_{E}},$$
$$\left( {1}/{n!}\; \right)U_{n}^{-1}V_{n}^{-1}{{\left| {{e}^{x}} \right|}^{-1}}{{x}^{p}}=\left[ p,\to  \right]{{\left( 1,\log \left( 1+x \right) \right)}_{E}},$$
$${{\left( 1,\log \left( 1+x \right) \right)}_{E}}{{\left( {{e}^{n\beta }},x \right)}_{E}}{{\left( 1,{{e}^{x}}-1 \right)}_{E}}={{\left( {{\left( 1+x \right)}^{n\beta }},x \right)}_{E}},$$
then
$${{V}_{n}}{{U}_{n}}{{E}^{n\beta }}U_{n}^{-1}V_{n}^{-1}={{\left( {{\left( 1+x \right)}^{n\beta }},x \right)}^{T}}{{I}_{n}}. \qquad\square $$

For example,
$${{G}_{3}}=\left( \begin{matrix}
   1 & 0 & 0 & 0  \\
   -3 & 1 & 0 & 0  \\
   3 & -2 & 1 & 0  \\
   -1 & 1 & -1 & 1  \\
\end{matrix} \right)\left( \begin{matrix}
   1 & 3 & 3 & 1  \\
   0 & 1 & 3 & 3  \\
   0 & 0 & 1 & 3  \\
   0 & 0 & 0 & 1  \\
\end{matrix} \right)\left( \begin{matrix}
   1 & 0 & 0 & 0  \\
   3 & 1 & 0 & 0  \\
   3 & 2 & 1 & 0  \\
   1 & 1 & 1 & 1  \\
\end{matrix} \right).$$
{\bfseries Theorem 6.3.} \emph{}
$$G_{n}^{\beta }{{x}^{p}}=\sum\limits_{m=0}^{n}{\left( \begin{matrix}
   -n\beta +p  \\
   m  \\
\end{matrix} \right)\left( \begin{matrix}
   n\beta +n-p  \\
   n-m  \\
\end{matrix} \right){{x}^{m}}}.$$
{\bfseries Proof.} We will use the factorial representation of binomial coefficients. Then
$${{\left( {{\left( 1+x \right)}^{n\beta }},x \right)}^{T}}{{V}_{n}}{{x}^{p}}=\sum\limits_{i=0}^{n}{\left( \begin{matrix}
   n\beta +n-p  \\
   n-i  \\
\end{matrix} \right){{x}^{i}}},$$
$$\left[ m,\to  \right]V_{n}^{-1}=\sum\limits_{i=0}^{m}{\left( \begin{matrix}
   m-n-1  \\
   m-i  \\
\end{matrix} \right){{x}^{i}}}=\sum\limits_{i=0}^{m}{{{\left( -1 \right)}^{m-i}}}\left( \begin{matrix}
   n-i  \\
   m-i  \\
\end{matrix} \right){{x}^{i}},$$
$$\left[ {{x}^{m}} \right]V_{n}^{-1}{{\left( {{\left( 1+x \right)}^{n\beta }},x \right)}^{T}}{{V}_{n}}{{x}^{p}}=$$
$$=\sum\limits_{i=0}^{m}{{{\left( -1 \right)}^{m-i}}}\left( \begin{matrix}
   n-i  \\
   m-i  \\
\end{matrix} \right)\left( \begin{matrix}
   n\beta +n-p  \\
   n-i  \\
\end{matrix} \right)\frac{\left( n\beta +m-p \right)!}{\left( n\beta +m-p \right)!}=$$
$$=\left( \begin{matrix}
   n\beta +n-p  \\
   n-m  \\
\end{matrix} \right)\sum\limits_{i=0}^{m}{{{\left( -1 \right)}^{m-i}}}\left( \begin{matrix}
   n\beta +m-p  \\
   m-i  \\
\end{matrix} \right)=$$
$$=\left( \begin{matrix}
   n\beta +n-p  \\
   n-m  \\
\end{matrix} \right){{\left( -1 \right)}^{m}}\left( \begin{matrix}
   n\beta +m-p-1  \\
   m  \\
\end{matrix} \right)=\left( \begin{matrix}
   n\beta +n-p  \\
   n-m  \\
\end{matrix} \right)\left( \begin{matrix}
   -n\beta +p  \\
   m  \\
\end{matrix} \right). \quad\square $$

We introduce the matrices ${{X}_{n}}=V_{n}^{-1}{{\left( x,x \right)}^{T}}{{V}_{n}}$. We find:
$${{X}_{n}}{{x}^{0}}=\frac{1-x-{{\left( 1-x \right)}^{n+1}}}{x} ,   \qquad{{X}_{n}}{{x}^{p}}={{x}^{p-1}}\left( 1-x \right).$$
Then
$$G_{n}^{\beta }={{\left( {{I}_{n}}+{{X}_{n}} \right)}^{n\beta }}=\sum\limits_{m=0}^{n}{\left( \begin{matrix}
   n\beta   \\
   m  \\
\end{matrix} \right)X_{n}^{m}}.$$
For example,
$${{G}_{3}}={{I}_{3}}+3\left( \begin{matrix}
   3 & 1 & 0 & 0  \\
   -6 & -1 & 1 & 0  \\
   4 & 0 & -1 & 1  \\
   -1 & 0 & 0 & -1  \\
\end{matrix} \right)+3\left( \begin{matrix}
   3 & 2 & 1 & 0  \\
   -8 & -5 & -2 & 1  \\
   7 & 4 & 1 & -2  \\
   -2 & -1 & 0 & 1  \\
\end{matrix} \right)+\left( \begin{matrix}
   1 & 1 & 1 & 1  \\
   -3 & -3 & -3 & -3  \\
   3 & 3 & 3 & 3  \\
   -1 & -1 & -1 & -1  \\
\end{matrix} \right),$$
$$G_{3}^{-1}={{I}_{3}}-3\left( \begin{matrix}
   3 & 1 & 0 & 0  \\
   -6 & -1 & 1 & 0  \\
   4 & 0 & -1 & 1  \\
   -1 & 0 & 0 & -1  \\
\end{matrix} \right)+6\left( \begin{matrix}
   3 & 2 & 1 & 0  \\
   -8 & -5 & -2 & 1  \\
   7 & 4 & 1 & -2  \\
   -2 & -1 & 0 & 1  \\
\end{matrix} \right)-10\left( \begin{matrix}
   1 & 1 & 1 & 1  \\
   -3 & -3 & -3 & -3  \\
   3 & 3 & 3 & 3  \\
   -1 & -1 & -1 & -1  \\
\end{matrix} \right).$$
Thus, ${{I}_{n}}+{{X}_{n}}=G_{n}^{{1}/{n}\;}={{U}_{n}}EU_{n}^{-1}$. For example,
$$G_{2}^{{1}/{2}\;}=\left( \begin{matrix}
   3 & 1 & 0  \\
   -3 & 0 & 1  \\
   1 & 0 & 0  \\
\end{matrix} \right),   \qquad G_{3}^{{1}/{3}\;}\left( \begin{matrix}
   4 & 1 & 0 & 0  \\
   -6 & 0 & 1 & 0  \\
   4 & 0 & 0 & 1  \\
   -1 & 0 & 0 & 0  \\
\end{matrix} \right),  \qquad G_{4}^{{1}/{4}\;}=\left( \begin{matrix}
   5 & 1 & 0 & 0 & 0  \\
   -10 & 0 & 1 & 0 & 0  \\
   10 & 0 & 0 & 1 & 0  \\
   -5 & 0 & 0 & 0 & 1  \\
   1 & 0 & 0 & 0 & 0  \\
\end{matrix} \right).$$

Note the identities  that follow from the Theorem 2.4. 
$$\left( {{\left( 1-x \right)}^{-m}},x \right)G_{n}^{\beta }\left( {{\left( 1-x \right)}^{m}},x \right){{I}_{n-m}}=G_{n-m}^{\frac{n\beta }{n-m}}.$$
$$\left( {{\left( 1-x \right)}^{-m}},x \right)G_{n}^{{1}/{n}\;}\left( {{\left( 1-x \right)}^{m}},x \right){{I}_{n-m}}=G_{n-m}^{{1}/{\left( n-m \right)}\;}.$$
{\bfseries Example 6.1.} Let $_{\left( \beta  \right)}a\left( x \right)$ is the generalized binomial series. Then
$$\left[ {{P}_{n}}\left( x \right) \right]\left( 1,xa\left( x \right) \right)={{x}^{n}},   \qquad\left[ {{P}_{n}}\left( x \right) \right]\left( 1,x{}_{\left( 1 \right)}a\left( x \right) \right)=x,$$
  $$\left[ {{P}_{n}}\left( x \right) \right]\left( a\left( x \right),xa\left( x \right) \right)={{x}^{n-1}}, \qquad\left[ {{P}_{n}}\left( x \right) \right]\left( _{\left( 1 \right)}a\left( x \right),x{}_{\left( 1 \right)}a\left( x \right) \right)=1.$$ 
Hence,
$$G_{n}^{\beta }{{x}^{n}}=\left[ {{P}_{n}}\left( x \right) \right]\left( 1+x{{\left( \log {}_{\left( \beta  \right)}{{a}^{\beta }}\left( x \right) \right)}^{\prime }},{{x}_{\left( \beta  \right)}}a\left( x \right) \right),$$
$$G_{n}^{\beta }x=\left[ {{P}_{n}}\left( x \right) \right]\left( 1+x{{\left( \log {}_{\left( \beta +1 \right)}{{a}^{\beta }}\left( x \right) \right)}^{\prime }},{{x}_{\left( \beta +1 \right)}}a\left( x \right) \right),$$
$$G_{n}^{\beta }{{x}^{n-1}}=\left[ {{P}_{n}}\left( x \right) \right]\left( _{\left( \beta  \right)}a\left( x \right)\left( 1+x{{\left( \log {}_{\left( \beta  \right)}{{a}^{\beta }}\left( x \right) \right)}^{\prime }} \right),{{x}_{\left( \beta  \right)}}a\left( x \right) \right),$$
$$G_{n}^{\beta }{{x}^{0}}=\left[ {{P}_{n}}\left( x \right) \right]\left( _{\left( \beta +1 \right)}a\left( x \right)\left( 1+x{{\left( \log {}_{\left( \beta +1 \right)}{{a}^{\beta }}\left( x \right) \right)}^{\prime }} \right),{{x}_{\left( \beta +1 \right)}}a\left( x \right) \right).$$
Since $G_{n}^{\beta }{{x}^{0}}=G_{n}^{\beta +1}{{x}^{n}}$, then the identity is manifested here:
$$_{\left( \beta  \right)}a\left( x \right)\left( 1+x{{\left( \log {}_{\left( \beta  \right)}{{a}^{\beta -1}}\left( x \right) \right)}^{\prime }} \right)=1+x{{\left( \log {}_{\left( \beta  \right)}{{a}^{\beta }}\left( x \right) \right)}^{\prime }}.$$
To match the identities $G_{n}^{-\beta }{{x}^{0}}={{J}_{n}}G_{n}^{\beta }{{x}^{n}}$, $G_{n}^{-\beta }x={{J}_{n}}G_{n}^{\beta }{{x}^{n-1}}$ with the definition of polynomials ${{\left( -1 \right)}^{n}}{{J}_{n}}{{g}_{n}}\left( x \right)$, we use the identity $_{\left( 1-\beta  \right)}a\left( x \right)={}_{\left( \beta  \right)}{{a}^{-1}}\left( -x \right)$.
\section{Polynomials $_{\left( \beta  \right)}{{h}_{n}}\left( x \right)$. }
Denote
$$\left[ {{P}_{n}}\left( x \right) \right]{{\left( b\left( x{}_{\left( \beta  \right)}{{a}^{\beta }}\left( x \right) \right)\left( 1+x{{\left( \log {}_{\left( \beta  \right)}{{a}^{\beta }}\left( x \right) \right)}^{\prime }} \right),x{}_{\left( \beta  \right)}a\left( x \right) \right)}_{E}}{{=}_{\left( \beta  \right)}}{{h}_{n}}\left( x \right).$$
We introduce the matrices $H_{n}^{\beta }={{F}_{n}}{{E}^{n\beta }}F_{n}^{-1}$. Then $H_{n}^{\beta }{{h}_{n}}\left( x \right)={}_{\left( \beta  \right)}{{h}_{n}}\left( x \right)$. For example,
$${{H}_{1}}=\frac{1}{2}\left( \begin{matrix}
   3 & 1  \\
   -1 & 1  \\
\end{matrix} \right), \quad{{H}_{2}}=\frac{1}{6}\left( \begin{matrix}
   15 & 5 & 1  \\
   -12 & 2 & 4  \\
   3 & -1 & 1  \\
\end{matrix} \right),  \quad{{H}_{3}}=\frac{1}{20}\left( \begin{matrix}
   84 & 28 & 7 & 1  \\
   -108 & -4 & 15 & 9  \\
   54 & -6 & -1 & 9  \\
   -10 & 2 & -1 & 1  \\
\end{matrix} \right).$$
Sum of elements of each column of the matrix $H_{n}^{\beta }$ is equal to one for the same reason as for the matrix $G_{n}^{\beta }$.
\\{\bfseries Theorem 7.1.} 
$$H_{n}^{-\beta }={{J}_{n}}H_{n}^{\beta }{{J}_{n}}.$$
{\bfseries Proof.} By the Theorem 3.1.
$${{J}_{n}}{{F}_{n}}{{E}^{n\beta }}F_{n}^{-1}{{J}_{n}}={{F}_{n}}{{E}^{n+1}}\left( 1,-x \right){{E}^{n\beta }}{{E}^{n+1}}\left( 1,-x \right)F_{n}^{-1}={{F}_{n}}{{E}^{-n\beta }}F_{n}^{-1}.\quad\square $$ 

Matrix $H_{n}^{\beta }$ can be represented as
$$H_{n}^{\beta }={{S}_{n}}G_{n}^{\beta }S_{n}^{-1}=V_{n}^{-1}{{C}_{n}}{{\left( {{\left( 1+x \right)}^{n\beta }},x \right)}^{T}}C_{n}^{-1}{{V}_{n}}.$$
 For example, 
$${{H}_{3}}=\left( \begin{matrix}
   1 & 0 & 0 & 0  \\
   -3 & 1 & 0 & 0  \\
   3 & -2 & 1 & 0  \\
   -1 & 1 & -1 & 1  \\
\end{matrix} \right)\left( \begin{matrix}
   1 & 0 & 0 & 0  \\
   0 & 4 & 0 & 0  \\
   0 & 0 & 10 & 0  \\
   0 & 0 & 0 & 20  \\
\end{matrix} \right)\left( \begin{matrix}
   1 & 3 & 3 & 1  \\
   0 & 1 & 3 & 3  \\
   0 & 0 & 1 & 3  \\
   0 & 0 & 0 & 1  \\
\end{matrix} \right)\left( \begin{matrix}
   1 & 0 & 0 & 0  \\
   0 & \frac{1}{4} & 0 & 0  \\
   0 & 0 & \frac{1}{10} & 0  \\
   0 & 0 & 0 & \frac{1}{20}  \\
\end{matrix} \right)\left( \begin{matrix}
   1 & 0 & 0 & 0  \\
   3 & 1 & 0 & 0  \\
   3 & 2 & 1 & 0  \\
   1 & 1 & 1 & 1  \\
\end{matrix} \right).$$

Denote
$${{t}_{n}}\left( \varphi |\beta ,x \right)=\sum\limits_{m=0}^{n}{\left( \begin{matrix}
   \varphi   \\
   m  \\
\end{matrix} \right)\left( \begin{matrix}
   \beta   \\
   n-m  \\
\end{matrix} \right){{x}^{m}}}.$$
{\bfseries Theorem7. 2.} \emph{}
$$H_{n}^{\beta }{{x}^{p}}=\sum\limits_{m=p}^{n}{{{\left( \begin{matrix}
   n+m  \\
   m  \\
\end{matrix} \right)}^{-1}}\left( \begin{matrix}
   n-p  \\
   n-m  \\
\end{matrix} \right)}{{\left( 1-x \right)}^{n-m}}{{t}_{m}}\left( -n\beta +n+m|n\beta ,x \right).$$
{\bfseries Proof.} If ${{c}_{p}}\left( x \right)$ is a polynomial of degree $p\le n$, then $V_{n}^{-1}{{c}_{p}}\left( x \right)={{\left( 1-x \right)}^{n-p}}V_{p}^{-1}{{c}_{p}}\left( x \right)$. Then
$$\frac{1}{n!}\left[ {{x}^{m}} \right]V_{p}^{-1}{{C}_{n}}{{\left( {{\left( 1+x \right)}^{n\beta }},x \right)}^{T}}{{x}^{p}}=$$
$$=\sum\limits_{i=0}^{m}{{{\left( -1 \right)}^{m-i}}}\left( \begin{matrix}
   p-i  \\
   m-i  \\
\end{matrix} \right)\left( \begin{matrix}
   n\beta   \\
   p-i  \\
\end{matrix} \right)\left( \begin{matrix}
   n+i  \\
   i  \\
\end{matrix} \right)\frac{\left( n\beta +m-p \right)!}{\left( n\beta +m-p \right)!}=$$
$$=\left( \begin{matrix}
   n\beta   \\
   p-m  \\
\end{matrix} \right)\sum\limits_{i=0}^{m}{{{\left( -1 \right)}^{m-i}}\left( \begin{matrix}
   n+i  \\
   i  \\
\end{matrix} \right)}\left( \begin{matrix}
   n\beta +m-p  \\
   m-i  \\
\end{matrix} \right)=$$
$$=\left( \begin{matrix}
   n\beta   \\
   p-m  \\
\end{matrix} \right){{\left( -1 \right)}^{m}}\left( \begin{matrix}
   n\beta +m-p-n-1  \\
   m  \\
\end{matrix} \right)=\left( \begin{matrix}
   n\beta   \\
   p-m  \\
\end{matrix} \right)\left( \begin{matrix}
   -n\beta +n+p  \\
   m  \\
\end{matrix} \right),$$
$$V_{n}^{-1}{{C}_{n}}{{\left( {{\left( 1+x \right)}^{n\beta }},x \right)}^{T}}{{x}^{p}}=n!{{\left( 1-x \right)}^{n-p}}{{t}_{p}}\left( -n\beta +n+p|n\beta ,x \right).$$
It remains to add that
$$C_{n}^{-1}{{V}_{n}}{{x}^{p}}=\frac{1}{n!}\sum\limits_{m=p}^{n}{{{\left( \begin{matrix}
   n+m  \\
   m  \\
\end{matrix} \right)}^{-1}}\left( \begin{matrix}
   n-p  \\
   n-m  \\
\end{matrix} \right)}{{x}^{m}}. \qquad\square $$

In particular,
$$H_{n}^{\beta }{{x}^{n}}={{\left( \begin{matrix}
   2n  \\
   n  \\
\end{matrix} \right)}^{-1}}\sum\limits_{m=0}^{n}{\left( \begin{matrix}
   -n\beta +2n  \\
   m  \\
\end{matrix} \right)\left( \begin{matrix}
   n\beta   \\
   n-m  \\
\end{matrix} \right){{x}^{m}}},$$
$$ H_{n}^{\beta }{{x}^{0}}={{\left( \begin{matrix}
   2n  \\
   n  \\
\end{matrix} \right)}^{-1}}\sum\limits_{m=0}^{n}{\left( \begin{matrix}
   -n\beta   \\
   m  \\
\end{matrix} \right)\left( \begin{matrix}
   n\beta +2n  \\
   n-m  \\
\end{matrix} \right){{x}^{m}}}.$$
{\bfseries Example 7.1.} Let $_{\left( \beta  \right)}a\left( x \right)$ is the generalized binomial series. Then, considering that ${{\left( 1,xa\left( x \right) \right)}^{-1}}=\left( 1,x{}_{\left( 2 \right)}a\left( -x \right) \right)$,
$$\left[ {{P}_{n}}\left( x \right) \right]{{\left( 1,xa\left( x \right) \right)}_{E}}={{S}_{n}}{{x}^{n}}=\frac{\left( 2n \right)!}{n!}{{x}^{n}}, \quad\left[ {{P}_{n}}\left( x \right) \right]{{\left( 1,x{}_{\left( 2 \right)}a\left( x \right) \right)}_{E}}=\frac{\left( 2n \right)!}{n!}x,$$
$$\left[ {{P}_{n}}\left( x \right) \right]{{\left( {{\left( xa\left( x \right) \right)}^{\prime }},xa\left( x \right) \right)}_{E}}=\frac{\left( 2n \right)!}{n!}{{x}^{n-1}},\quad
\left[ {{P}_{n}}\left( x \right) \right]{{\left( {{\left( {{x}_{\left( 2 \right)}}a\left( x \right) \right)}^{\prime }},{{x}_{\left( 2 \right)}}a\left( x \right) \right)}_{E}}=\frac{\left( 2n \right)!}{n!},$$
$$\left[ {{P}_{n}}\left( x \right) \right]{{\left( a\left( x \right),xa\left( x \right) \right)}_{E}}={{S}_{n}}{{x}^{n-1}}=\frac{\left( 2n \right)!}{n!2}\left( 1+x \right){{x}^{n-1}},$$
$$\left[ {{P}_{n}}\left( x \right) \right]{{\left( 1+x{{\left( \log {}_{\left( 2 \right)}a\left( x \right) \right)}^{\prime }},x{}_{\left( 2 \right)}a\left( x \right) \right)}_{E}}=\frac{\left( 2n \right)!}{n!2}\left( 1+x \right).$$
Hence,
$$\frac{\left( 2n \right)!}{n!}H_{n}^{\beta }{{x}^{n}}=\left[ {{P}_{n}}\left( x \right) \right]{{\left( 1+x{{\left( \log {}_{\left( \beta  \right)}{{a}^{\beta }}\left( x \right) \right)}^{\prime }},{{x}_{\left( \beta  \right)}}a\left( x \right) \right)}_{E}},$$
$$\frac{\left( 2n \right)!}{n!}H_{n}^{\beta }{{x}^{0}}=\frac{\left( 2n \right)!}{n!}H_{n}^{\beta +2}{{x}^{n}}=\left[ {{P}_{n}}\left( x \right) \right]{{\left( 1+x{{\left( \log {}_{\left( \beta +2 \right)}{{a}^{\beta +2}}\left( x \right) \right)}^{\prime }},{{x}_{\left( \beta +2 \right)}}a\left( x \right) \right)}_{E}},$$
$$\frac{\left( 2n \right)!}{n!}H_{n}^{\beta }x=\left[ {{P}_{n}}\left( x \right) \right]{{\left( 1+x{{\left( \log {}_{\left( \beta +2 \right)}{{a}^{\beta }}\left( x \right) \right)}^{\prime }},{{x}_{\left( \beta +2 \right)}}a\left( x \right) \right)}_{E}},$$
$$\frac{\left( 2n \right)!}{n!}H_{n}^{\beta }{{x}^{n-1}}=\left[ {{P}_{n}}\left( x \right) \right]{{\left( _{\left( \beta  \right)}a\left( x \right)\left( 1+x{{\left( \log {}_{\left( \beta  \right)}{{a}^{\beta +1}}\left( x \right) \right)}^{\prime }} \right),{{x}_{\left( \beta  \right)}}a\left( x \right) \right)}_{E}},$$
$$\frac{\left( 2n \right)!}{n!2}H_{n}^{\beta }\left( 1+x \right){{x}^{n-1}}=\left[ {{P}_{n}}\left( x \right) \right]{{\left( _{\left( \beta  \right)}a\left( x \right)\left( 1+x{{\left( \log {}_{\left( \beta  \right)}{{a}^{\beta }}\left( x \right) \right)}^{\prime }} \right),{{x}_{\left( \beta  \right)}}a\left( x \right) \right)}_{E}},$$
$$\frac{\left( 2n \right)!}{n!2}H_{n}^{\beta }\left( 1+x \right)=\left[ {{P}_{n}}\left( x \right) \right]{{\left( 1+x{{\left( \log {}_{\left( \beta +2 \right)}{{a}^{\beta +1}}\left( x \right) \right)}^{\prime }},x{}_{\left( \beta +2 \right)}a\left( x \right) \right)}_{E}}.$$
Note that
$$\left( \begin{matrix}
   2n-1  \\
   n-1  \\
\end{matrix} \right)H_{n}^{\beta }\left( 1+x \right){{x}^{n-1}}=\left( 1-x \right)\sum\limits_{m=0}^{n-1}{\left( \begin{matrix}
   -n\beta +2n-1  \\
   m  \\
\end{matrix} \right)\left( \begin{matrix}
   n\beta   \\
   n-1-m  \\
\end{matrix} \right){{x}^{m}}}+$$
$$+\sum\limits_{m=0}^{n}{\left( \begin{matrix}
   -n\beta +2n  \\
   m  \\
\end{matrix} \right)\left( \begin{matrix}
   n\beta   \\
   n-m  \\
\end{matrix} \right){{x}^{m}}}=\sum\limits_{m=0}^{n}{\left( \begin{matrix}
   -n\beta +2n-1  \\
   m  \\
\end{matrix} \right)\left( \begin{matrix}
   n\beta +1  \\
   n-m  \\
\end{matrix} \right){{x}^{m}}},$$
$$\left( \begin{matrix}
   2n-1  \\
   n-1  \\
\end{matrix} \right)H_{n}^{\beta }\left( 1+x \right)=\sum\limits_{m=0}^{n}{\left( \begin{matrix}
   -n\beta +1  \\
   m  \\
\end{matrix} \right)\left( \begin{matrix}
   n\beta +2n-1  \\
   n-m  \\
\end{matrix} \right){{x}^{m}}}.$$
To match the identities $H_{n}^{-\beta }{{x}^{0}}={{J}_{n}}H_{n}^{\beta }{{x}^{n}}$, $H_{n}^{-\beta }x={{J}_{n}}H_{n}^{\beta }{{x}^{n-1}}$ with the definition of polynomials ${{\left( -1 \right)}^{n}}{{J}_{n}}{{h}_{n}}\left( x \right)$, we use the identities  ${{\left( 1,x{}_{\left( \beta  \right)}a\left( x \right) \right)}^{-1}}=\left( 1,{{x}_{\left( \beta -1 \right)}}{{a}^{-1}}\left( x \right) \right)$,  ${}_{\left( \beta -1 \right)}{{a}^{-1}}\left( -x \right){{=}_{\left( 2-\beta  \right)}}a\left( x \right)$.
\section{ Matrices ${{\tilde{U}}_{n}}$, ${{\tilde{F}}_{n}}$}

Let, as before, $\left[ {{P}_{n}}\left( x \right) \right]\left( 1,xa\left( x \right) \right)={{\alpha }_{n}}\left( x \right)$, $\left[ n,\to  \right]{{\left( 1,\log a\left( x \right) \right)}_{E}}={{u}_{n}}\left( x \right)$. For the sequence of polynomials ${{c}_{n}}\left( x \right)$ such that ${{c}_{0}}\left( 0 \right)=1$, ${{c}_{n}}\left( 0 \right)=0$ we denote: ${{\tilde{c}}_{0}}\left( x \right)=1$, ${{\tilde{c}}_{n}}\left( x \right)=\left( {1}/{x}\; \right){{c}_{n}}\left( x \right)$. We introduce the matrices ${{\tilde{U}}_{n}}$, $n>0$:
$${{\tilde{U}}_{n}}={{\left( x,x \right)}^{T}}{{U}_{n}}\left( x,x \right),   \qquad\tilde{U}_{n}^{-1}={{\left( x,x \right)}^{T}}U_{n}^{-1}\left( x,x \right),$$
$${{\tilde{U}}_{n}}{{x}^{p}}=\frac{1}{n!}{{\left( 1-x \right)}^{n-1-p}}{{\tilde{A}}_{p+1}}\left( x \right),  \quad\tilde{U}_{n}^{-1}{{x}^{p}}={{\left( x-1 \right)}_{p}}{{\left[ x+1 \right]}_{n-p-1}}, \quad 0\le p<n.$$
For example, 
$$\tilde{U}_{4}^{{}}=\frac{1}{4!}\left( \begin{matrix}
   1 & 1 & 1 & 1  \\
   -3 & -1 & 3 & 11  \\
   3 & -1 & -3 & 11  \\
   -1 & 1 & -1 & 1  \\
\end{matrix} \right),  \qquad\tilde{U}_{4}^{-1}=\left( \begin{matrix}
   6 & -2 & 2 & -6  \\
   11 & -1 & -1 & 11  \\
   6 & 2 & -2 & -6  \\
   1 & 1 & 1 & 1  \\
\end{matrix} \right).$$
Then ${{\tilde{U}}_{n}}{{\tilde{u}}_{n}}\left( x \right)={{\tilde{\alpha }}_{n}}\left( x \right)$. Since
$$\frac{{{{\tilde{\alpha }}}_{n}}\left( x \right)}{{{\left( 1-x \right)}^{n+1}}}=\sum\limits_{m=0}^{\infty }{\frac{{{u}_{n}}\left( m+1 \right)}{n!}}{{x}^{m}},  \qquad{{\tilde{\alpha }}_{n}}\left( x \right)={{U}_{n}}E{{u}_{n}}\left( x \right),$$ 
then the matrices ${{\tilde{U}}_{n}}$, $\tilde{U}_{n}^{-1}$ can be represented as
 $${{\tilde{U}}_{n}}={{U}_{n}}E\left( x,x \right){{I}_{n-1}},   \qquad\tilde{U}_{n}^{-1}={{\left( x,x \right)}^{T}}{{E}^{-1}}U_{n}^{-1}{{I}_{n-1}}.$$

 We introduce the matrices ${{\tilde{J}}_{n}}={{J}_{n-1}}$, ${{\tilde{I}}_{n}}={{I}_{n-1}}$, ${{\tilde{V}}_{n}}={{V}_{n-1}}$. Then
$$\tilde{U}_{n}^{-1}\tilde{V}_{n}^{-1}{{x}^{p}}=\frac{n!}{\left( p+1 \right)!}\sum\limits_{m=0}^{p}{s\left( p+1,\text{ }m+1 \right){{x}^{m}}},
\quad{{\tilde{V}}_{n}}{{\tilde{U}}_{n}}{{x}^{p}}=\frac{1}{n!}\sum\limits_{m=0}^{p}{m!S\left( p+1,\text{ }m+1 \right)}\text{ }{{x}^{m}}.$$
For example, $\tilde{U}_{4}^{-1}\tilde{V}_{4}^{-1}$, ${{\tilde{V}}_{4}}{{\tilde{U}}_{4}}$:
$$4!\left( \begin{matrix}
   1 & -1 & \text{  }2 & -6  \\
   0 & \text{  }1 & -3 & \text{ }11  \\
   0 & \text{  }0 & \text{  }1 & -6  \\
   0 & \text{  }0 & \text{  }0 & \text{  }1  \\
\end{matrix} \right)\left( \begin{matrix}
   1 & 0 & 0 & 0  \\
   0 & \frac{1}{2} & 0 & 0  \\
   0 & 0 & \frac{1}{3!} & 0  \\
   0 & 0 & 0 & \frac{1}{4!}  \\
\end{matrix} \right), \qquad\frac{1}{4!}\left( \begin{matrix}
   1 & 0 & 0 & 0  \\
   0 & 2 & 0 & 0  \\
   0 & 0 & 3! & 0  \\
   0 & 0 & 0 & 4!  \\
\end{matrix} \right)\left( \begin{matrix}
   1 & 1 & 1 & 1  \\
   0 & 1 & 3 & 7  \\
   0 & 0 & 1 & 6  \\
   0 & 0 & 0 & 1  \\
\end{matrix} \right).$$
{\bfseries Theorem 8.1.} \emph{}
$${{\tilde{U}}_{n}}\left( 1,-x \right)\tilde{U}_{n}^{-1}={{\left( -1 \right)}^{n-1}}{{\tilde{J}}_{n}}.$$
{\bfseries Proof.}
$$\left( 1,-x \right){{\left( x-1 \right)}_{p}}{{\left[ x+1 \right]}_{n-p-1}}={{\left( -x-1 \right)}_{p}}{{\left[ -x+1 \right]}_{n-p-1}}={{\left( -1 \right)}^{n-1}}{{\left( x-1 \right)}_{n-p-1}}{{\left[ x+1 \right]}_{p}}.\,\square $$
{\bfseries Example 8.1.} This example illustrates the advantage of the matrices ${{\tilde{U}}_{n}}$ over the matrices ${{U}_{n}}$ in certain cases. Denote
$$\left[ {{P}_{n}}\left( x \right) \right]\left( 1,x{{a}^{m}}\left( x \right) \right)=\alpha _{n}^{\left( m \right)}\left( x \right), \qquad{{W}_{\left( n,m \right)}}={{\tilde{U}}_{n}}\left( m,mx \right)\tilde{U}_{n}^{-1}.$$
Then ${{W}_{\left( n,m \right)}}{{\tilde{\alpha }}_{n}}\left( x \right)=\tilde{\alpha }_{n}^{\left( m \right)}\left( x \right)$. We build the matrix ${{\left( a\left( x \right),x \right)}_{m}}$ by the rule $\left[ n,\to  \right]{{\left( a\left( x \right),x \right)}_{m}}=\left[ mn+m-1,\to  \right]\left( a\left( x \right),x \right)$. For example, 
$${{\left( a\left( x \right),x \right)}_{2}}=\left( \begin{matrix}
   {{a}_{1}} & {{a}_{0}} & 0 & 0 & \cdots   \\
   {{a}_{3}} & {{a}_{2}} & {{a}_{1}} & {{a}_{0}} & \cdots   \\
   {{a}_{5}} & {{a}_{4}} & {{a}_{3}} & {{a}_{2}} & \cdots   \\
   {{a}_{7}} & {{a}_{6}} & {{a}_{5}} & {{a}_{4}} & \cdots   \\
   \vdots  & \vdots  & \vdots  & \vdots  & \ddots   \\
\end{matrix} \right),  \qquad{{\left( a\left( x \right),x \right)}_{3}}=\left( \begin{matrix}
   {{a}_{2}} & {{a}_{1}} & {{a}_{0}} & 0 & \cdots   \\
   {{a}_{5}} & {{a}_{4}} & {{a}_{3}} & {{a}_{2}} & \cdots   \\
   {{a}_{8}} & {{a}_{7}} & {{a}_{6}} & {{a}_{5}} & \cdots   \\
   {{a}_{11}} & {{a}_{10}} & {{a}_{9}} & {{a}_{8}} & \cdots   \\
   \vdots  & \vdots  & \vdots  & \vdots  & \ddots   \\
\end{matrix} \right).$$
{\bfseries Theorem 8.2.} 
$${{W}_{\left( n,m \right)}}={{\left( w_{m}^{n+1}\left( x \right),x \right)}_{m}}{{\tilde{I}}_{n}} ,  \qquad w_{m}^{n+1}\left( x \right)={{\left( \frac{1-{{x}^{m}}}{1-x} \right)}^{n+1}}.$$
{\bfseries Proof.} Since
$$\left[ {{x}^{p}} \right]\frac{\tilde{\alpha }_{n}^{\left( m \right)}\left( x \right)}{{{\left( 1-x \right)}^{n+1}}}=\left[ {{x}^{mp+m-1}} \right]\frac{{{{\tilde{\alpha }}}_{n}}\left( x \right)}{{{\left( 1-x \right)}^{n+1}}},
\quad\frac{{{{\tilde{\alpha }}}_{n}}\left( x \right)}{{{\left( 1-x \right)}^{n+1}}}=\frac{w_{m}^{n+1}\left( x \right){{{\tilde{\alpha }}}_{n}}\left( x \right)}{{{\left( 1-{{x}^{m}} \right)}^{n+1}}}=\sum\limits_{r=0}^{m-1}{\frac{{{x}^{r}}{{c}_{r}}\left( x \right)}{{{\left( 1-{{x}^{m}} \right)}^{n+1}}}},$$
$${{c}_{r}}\left( x \right)=\sum\limits_{p=0}^{\infty }{\left( \left[ {{x}^{mp+r}} \right]w_{m}^{n+1}\left( x \right){{{\tilde{\alpha }}}_{n}}\left( x \right) \right){{x}^{mp}}},$$ 
then
$$\frac{\tilde{\alpha }_{n}^{\left( m \right)}\left( {{x}^{m}} \right)}{{{\left( 1-{{x}^{m}} \right)}^{n+1}}}=\frac{{{c}_{m-1}}\left( x \right)}{{{\left( 1-{{x}^{m}} \right)}^{n+1}}}, \qquad\left[ {{x}^{p}} \right]\tilde{\alpha }_{n}^{\left( m \right)}\left( x \right)=\left[ {{x}^{mp+m-1}} \right]w_{m}^{n+1}\left( x \right){{\tilde{\alpha }}_{n}}\left( x \right),$$
$$\tilde{\alpha }_{n}^{\left( m \right)}\left( x \right)={{\left( w_{m}^{n+1}\left( x \right),x \right)}_{m}}{{\tilde{\alpha }}_{n}}\left( x \right).\qquad\square $$

For example, ${{W}_{(1,\text{2})}}=\left( 2 \right)$,  ${{W}_{(1,3)}}=\left( 3 \right)$,  ${{W}_{(1,4)}}=\left( 4 \right)$,
    $${{W}_{(2,2)}}=\left( \begin{matrix}
   3 & 1  \\
   1 & 3  \\
\end{matrix} \right),       \qquad{{W}_{(3,2)}}=\left( \begin{matrix}
   4 & 1 & 0  \\
   4 & 6 & 4  \\
   0 & 1 & 4  \\
\end{matrix} \right),  \qquad{{W}_{(4,2)}}=\left( \begin{matrix}
   5 & 1 & 0 & 0  \\
   10 & 10 & 5 & 1  \\
   1 & 5 & 10 & 10  \\
   0 & 0 & 1 & 5  \\
\end{matrix} \right),$$
   $${{W}_{(2,3)}}=\left( \begin{matrix}
   6 & 3  \\
   3 & 6  \\
\end{matrix} \right),      \qquad{{W}_{(3,3)}}=\left( \begin{matrix}
   10 & 4 & 1  \\
   16 & 19 & 16  \\
   1 & 4 & 10  \\
\end{matrix} \right),  \qquad{{W}_{(4,3)}}=\left( \begin{matrix}
   15 & 5 & 1 & 0  \\
   51 & 45 & 30 & 15  \\
   15 & 30 & 45 & 51  \\
   0 & 1 & 5 & 15  \\
\end{matrix} \right),$$
$${{W}_{(2,4)}}=\left( \begin{matrix}
   10 & 6  \\
   6 & 10  \\
\end{matrix} \right), \quad{{W}_{(3,4)}}=\left( \begin{matrix}
   20 & 10 & 4  \\
   40 & 44 & 40  \\
   4 & 10 & 20  \\
\end{matrix} \right), \quad{{W}_{\left( 4,4 \right)}}=\left( \begin{matrix}
   35 & 15 & 5 & 1  \\
   155 & 135 & 101 & 65  \\
   65 & 101 & 135 & 155  \\
   1 & 5 & 15 & 35  \\
\end{matrix} \right).$$
Since $\alpha _{n}^{\left( m \right)}\left( 1 \right)={{m}^{n}}{{\alpha }_{n}}\left( 1 \right)$, then the sum of  elements of each column of the matrix ${{W}_{\left( n,m \right)}}$ is equal to ${{m}^{n}}$. Note the identities
$${{W}_{\left( n,m \right)}}{{\tilde{A}}_{n}}\left( x \right)={{m}^{n}}{{\tilde{A}}_{n}}\left( x \right), \qquad{{W}_{\left( n,m \right)}}{{\tilde{J}}_{n}}={{\tilde{J}}_{n}}{{W}_{\left( n,\text{ }m \right)}},  \qquad{{W}_{\left( n,m \right)}}{{W}_{\left( n,p \right)}}={{W}_{\left( n,mp \right)}},$$
$$\left( {{\left( 1-x \right)}^{-p}},x \right){{W}_{\left( n,m \right)}}\left( {{\left( 1-x \right)}^{p}},x \right){{\tilde{I}}_{n-p}}={{W}_{\left( n-p,m \right)}}.$$
For example, 
$$\left( \begin{matrix}
   4 & 1 & 0  \\
   4 & 6 & 4  \\
   0 & 1 & 4  \\
\end{matrix} \right)\left( \begin{matrix}
   1  \\
   4  \\
   1  \\
\end{matrix} \right)=8\left( \begin{matrix}
   1  \\
   4  \\
   1  \\
\end{matrix} \right),  \quad\left( \begin{matrix}
   10 & 4 & 1  \\
   16 & 19 & 16  \\
   1 & 4 & 10  \\
\end{matrix} \right)\left( \begin{matrix}
   1  \\
   4  \\
   1  \\
\end{matrix} \right)=27\left( \begin{matrix}
   1  \\
   4  \\
   1  \\
\end{matrix} \right),$$
$$\left( \begin{matrix}
   4 & 1 & 0  \\
   4 & 6 & 4  \\
   0 & 1 & 4  \\
\end{matrix} \right)\left( \begin{matrix}
   4 & 1 & 0  \\
   4 & 6 & 4  \\
   0 & 1 & 4  \\
\end{matrix} \right)=\left( \begin{matrix}
   20 & 10 & 4  \\
   40 & 44 & 40  \\
   4 & 10 & 20  \\
\end{matrix} \right),$$
$$\left( \begin{matrix}
   1 & 0 & 0  \\
   1 & 1 & 0  \\
   1 & 1 & 1  \\
\end{matrix} \right)\left( \begin{matrix}
   4 & 1 & 0  \\
   4 & 6 & 4  \\
   0 & 1 & 4  \\
\end{matrix} \right)\left( \begin{matrix}
   \text{  }1 & \text{ }0  \\
   -1 & \text{  }1  \\
   \text{  }0 & -1  \\
\end{matrix} \right)=\left( \begin{matrix}
   3 & 1  \\
   1 & 3  \\
\end{matrix} \right),$$
$$\left( \begin{matrix}
   1 & 0 & 0  \\
   2 & 1 & 0  \\
   3 & 2 & 1  \\
\end{matrix} \right)\left( \begin{matrix}
   4 & 1 & 0  \\
   4 & 6 & 4  \\
   0 & 1 & 4  \\
\end{matrix} \right)\left( \begin{matrix}
   \text{  }1  \\
   -2  \\
   \text{  }1  \\
\end{matrix} \right)=\left( \begin{matrix}
   1 & 0  \\
   1 & 1  \\
\end{matrix} \right)\left( \begin{matrix}
   3 & 1  \\
   1 & 3  \\
\end{matrix} \right)\left( \begin{matrix}
   \text{  }1  \\
   -1  \\
\end{matrix} \right)=\left( 2 \right).$$
Since $\left( 1,a\left( x \right)-1 \right)\left( 1,{{\left( 1+x \right)}^{m}}-1 \right)=\left( 1,{{a}^{m}}\left( x \right)-1 \right)$, then the matrix ${{W}_{\left( n,m \right)}}$ can be represented as
$${{W}_{\left( n,m \right)}}=\tilde{V}_{n}^{-1}{{\left( \frac{{{\left( 1+x \right)}^{m}}-1}{x},{{\left( 1+x \right)}^{m}}-1 \right)}^{T}}{{\tilde{V}}_{n}}.$$
For example,
$$\left( \begin{matrix}
   4 & 1 & 0  \\
   4 & 6 & 4  \\
   0 & 1 & 4  \\
\end{matrix} \right)=\frac{1}{3!}\left( \begin{matrix}
   1 & 1 & 1  \\
   -2 & 0 & 4  \\
   1 & -1 & 1  \\
\end{matrix} \right)\left( \begin{matrix}
   2 & 0 & 0  \\
   0 & 4 & 0  \\
   0 & 0 & 8  \\
\end{matrix} \right)\left( \begin{matrix}
   2 & -1 & 2  \\
   3 & 0 & -3  \\
   1 & 1 & 1  \\
\end{matrix} \right)=$$
$$=\left( \begin{matrix}
   1 & 0 & 0  \\
   -2 & 1 & 0  \\
   1 & -1 & 1  \\
\end{matrix} \right)\left( \begin{matrix}
   2 & 1 & 0  \\
   0 & 4 & 4  \\
   0 & 0 & 8  \\
\end{matrix} \right)\left( \begin{matrix}
   1 & 0 & 0  \\
   2 & 1 & 0  \\
   1 & 1 & 1  \\
\end{matrix} \right).$$
Matrices $\left( {1}/{{{m}^{n}}}\; \right){{\left( {{W}_{\left( n,m \right)}} \right)}^{T}}$are known as amazing matrices [27, p.156]. They find application in various fields of mathematics [28] – [31].

Denote $\left[ {{P}_{n}}\left( x \right) \right]{{\left( 1,xa\left( x \right) \right)}_{E}}={{\varphi }_{n}}\left( x \right)$. We introduce the matrices ${{\tilde{F}}_{n}}$:
$${{\tilde{F}}_{n}}={{\left( x,x \right)}^{T}}{{F}_{n}}\left( x,x \right),   \qquad\tilde{F}_{n}^{-1}={{\left( x,x \right)}^{T}}U_{n}^{-1}\left( x,x \right),$$
$${{\tilde{F}}_{n}}{{x}^{p}}={{\left( 1-x \right)}^{2n+1}}\sum\limits_{m=0}^{\infty }{{{\left( m+1 \right)}^{p+1}}}\left( \begin{matrix}
   m+n+1  \\
   n  \\
\end{matrix} \right){{x}^{m}},$$
$$\tilde{F}_{n}^{-1}{{x}^{p}}=\frac{n!}{\left( 2n \right)!}{{\left( x-1 \right)}_{p}}{{\left[ x+n+1 \right]}_{n-p-1}}, \qquad 0\le p<n.$$
For example, 
$${{\tilde{F}}_{4}}=5\left( \begin{matrix}
   1 & 1 & 1 & 1  \\
   -3 & 3 & 15 & 39  \\
   3 & -9 & 9 & 171  \\
   -1 & 5 & -25 & 125  \\
\end{matrix} \right),  \qquad\tilde{F}_{4}^{-1}=\frac{4!}{8!}\left( \begin{matrix}
   210 & -30 & 10 & -6  \\
   107 & 19 & -13 & 11  \\
   18 & 10 & 2 & -6  \\
   1 & 1 & 1 & 1  \\
\end{matrix} \right).$$
Then ${{\tilde{F}}_{n}}{{\tilde{u}}_{n}}\left( x \right)={{\tilde{\varphi }}_{n}}\left( x \right)$. Since  (see Example 3.1.)
$$\frac{{{{\tilde{\varphi }}}_{n}}\left( x \right)}{{{\left( 1-x \right)}^{2n+1}}}=\sum\limits_{m=0}^{\infty }{\left( \begin{matrix}
   n+m  \\
   m  \\
\end{matrix} \right)\left( m+n+1 \right){{{\tilde{u}}}_{n}}\left( m+1 \right)}{{x}^{m}},  
\quad{{\tilde{\varphi }}_{n}}\left( x \right)={{F}_{n}}E\left( x+n \right){{\tilde{u}}_{n}}\left( x \right),$$
then the matrices ${{\tilde{F}}_{n}}$, $\tilde{F}_{n}^{-1}$ can be represented as
$${{\tilde{F}}_{n}}={{F}_{n}}E\left( x+n,x \right){{I}_{n-1}}={{F}_{n}}\left( x+n+1,x \right)E{{I}_{n-1}}, \quad\tilde{F}_{n}^{-1}={{\left( x+n,x \right)}^{-1}}{{E}^{-1}}F_{n}^{-1}{{I}_{n-1}}.$$
Note that since
$$\left[ {{x}^{n}} \right]{{F}_{n}}{{x}^{p}}\left( x+n+1 \right)=\left[ {{x}^{n}} \right]{{F}_{n}}E{{x}^{p}}\left( x+n \right)=0, \qquad p<n;$$
$${{F}_{n}}{{x}^{p}}={{\left( 1-x \right)}^{2n+1}}\sum\limits_{m=0}^{\infty }{{{m}^{p}}}\left( \begin{matrix}
   m+n  \\
   n  \\
\end{matrix} \right){{x}^{m}},\quad
{{F}_{n}}E{{I}_{n}}{{x}^{p}}={{\left( 1-x \right)}^{2n+1}}\sum\limits_{m=0}^{\infty }{{{\left( m+1 \right)}^{p}}}\left( \begin{matrix}
   m+n  \\
   n  \\
\end{matrix} \right){{x}^{m}},$$
then the identities for the $n$th elements of  columns of the matrices ${{F}_{n}}$, ${{F}_{n}}E{{I}_{n}}$ are manifested here:
$$\sum\limits_{m=0}^{n}{{{\left( -1 \right)}^{n-m}}\left( \begin{matrix}
   2n+1  \\
   n-m  \\
\end{matrix} \right){{m}^{p}}}\left( \begin{matrix}
   m+n  \\
   n  \\
\end{matrix} \right)={{\left( -1 \right)}^{n+p}}{{\left( n+1 \right)}^{p}},$$
$$\sum\limits_{m=0}^{n}{{{\left( -1 \right)}^{n-m}}}\left( \begin{matrix}
   2n+1  \\
   n-m  \\
\end{matrix} \right){{\left( m+1 \right)}^{p}}\left( \begin{matrix}
   m+n  \\
   n  \\
\end{matrix} \right)={{\left( -1 \right)}^{n+p}}{{n}^{p}},  \qquad p\le n.$$
{\bfseries Theorem 8.3.} 
$${{\tilde{F}}_{n}}{{E}^{n}}\left( 1,-x \right)\tilde{F}_{n}^{-1}={{\left( -1 \right)}^{n-1}}{{\tilde{J}}_{n}}.$$
{\bfseries Proof.} 
$${{E}^{n}}\left( 1,-x \right){{\left( x-1 \right)}_{p}}{{\left[ x+n+1 \right]}_{n-p-1}}={{\left( -x-n-1 \right)}_{p}}{{\left[ -x+1 \right]}_{n-p-1}}=$$
$$={{\left( -1 \right)}^{n-1}}{{\left( x-1 \right)}_{n-p-1}}{{\left[ x+n+1 \right]}_{p}}. \qquad\square $$

We introduce the matrices ${{\tilde{S}}_{n}}={{\tilde{V}}_{n}}{{\tilde{C}}_{n}}{{\tilde{V}}_{n}}$, ${{\tilde{C}}_{n}}{{x}^{p}}=\left( {\left( n+p+1 \right)!}/{\left( p+1 \right)!}\; \right){{x}^{p}}$. Then ${{\tilde{S}}_{n}}{{\tilde{\alpha }}_{n}}\left( x \right)={{\tilde{\varphi }}_{n}}\left( x \right)$.
\section{Polynomials $_{\left( \beta  \right)}{{\alpha }_{n}}\left( x \right)$, $_{\left( \beta  \right)}{{\varphi }_{n}}\left( x \right)$ }
Denote $\left[ {{P}_{n}}\left( x \right) \right]\left( 1,x{}_{\left( \beta  \right)}a\left( x \right) \right){{=}_{\left( \beta  \right)}}{{\alpha }_{n}}\left( x \right)$. We introduce the matrices $A_{n}^{\beta }={{\tilde{U}}_{n}}{{E}^{n\beta }}\tilde{U}_{n}^{-1}$. Then $A_{n}^{\beta }{{\tilde{\alpha }}_{n}}\left( x \right)={}_{\left( \beta  \right)}{{\tilde{\alpha }}_{n}}\left( x \right)$. For example, 
$${{A}_{2}}=\left( \begin{matrix}
   \text{  }2 & 1  \\
   -1 & 0  \\
\end{matrix} \right),   \quad{{A}_{3}}=\left( \begin{matrix}
   5 & {5}/{2}\; & 1  \\
   -6 & -2 & 0  \\
   2 & {1}/{2}\; & 0  \\
\end{matrix} \right),   \quad{{A}_{4}}=\left( \begin{matrix}
   14 & 7 & 3 & 1  \\
   -28 & {-35}/{3}\; & {-10}/{3}\; & 0  \\
   20 & {22}/{3}\; & {5}/{3}\; & 0  \\
   -5 & {-5}/{3}\; & {-1}/{3}\; & 0  \\
\end{matrix} \right).$$
{\bfseries Theorem 9.1.} 
$$A_{n}^{-\beta }={{\tilde{J}}_{n}}A_{n}^{\beta }{{\tilde{J}}_{n}}.$$
{\bfseries Proof.} By the Theorem 8.1.
$${{\tilde{J}}_{n}}{{\tilde{U}}_{n}}{{E}^{n\beta }}\tilde{U}_{n}^{-1}{{\tilde{J}}_{n}}={{\tilde{U}}_{n}}\left( 1,-x \right){{E}^{n\beta }}\left( 1,-x \right)\tilde{U}_{n}^{-1}={{\tilde{U}}_{n}}{{E}^{-n\beta }}\tilde{U}_{n}^{-1}. \qquad\square $$

Note the identity 
$$\left( {{\left( 1-x \right)}^{-m}},x \right)A_{n}^{\beta }\left( {{\left( 1-x \right)}^{m}},x \right){{\tilde{I}}_{n-m}}=A_{n-m}^{\frac{n\beta }{n-m}}.$$

Let $D$ is the matrix of derivation operator. We introduce the diagonal matrix \linebreak $\tilde{D}=D\left( x,x \right)$, $\tilde{D}{{x}^{n}}=\left( n+1 \right){{x}^{n}}$.
\\{\bfseries Theorem 9.2.} 
$$A_{n}^{\beta }=\tilde{V}_{n}^{-1}\tilde{D}{{\left( {{\left( 1+x \right)}^{n\beta }},x \right)}^{T}}{{\tilde{D}}^{-1}}{{\tilde{V}}_{n}}.$$
{\bfseries Proof.} Since
$$n!\left| {{e}^{x}} \right|{{\tilde{D}}^{-1}}{{\tilde{V}}_{n}}{{\tilde{U}}_{n}}{{x}^{p}}=\left[ p,\to  \right]{{\left( {{e}^{x}},{{e}^{x}}-1 \right)}_{E}},$$
$$\left( {1}/{n!}\; \right)\tilde{U}_{n}^{-1}\tilde{V}_{n}^{-1}{{\left| {{e}^{x}} \right|}^{-1}}\tilde{D}{{x}^{p}}=\left[ p,\to  \right]{{\left( {{\left( 1+x \right)}^{-1}},\log \left( 1+x \right) \right)}_{E}},$$
$${{\left( {{\left( 1+x \right)}^{-1}},\log \left( 1+x \right) \right)}_{E}}{{\left( {{e}^{n\beta }},x \right)}_{E}}{{\left( {{e}^{x}},{{e}^{x}}-1 \right)}_{E}}={{\left( {{\left( 1+x \right)}^{n\beta }},x \right)}_{E}},$$
then
$${{\tilde{V}}_{n}}{{\tilde{U}}_{n}}{{E}^{n\beta }}\tilde{U}_{n}^{-1}\tilde{V}_{n}^{-1}=\tilde{D}{{\left( {{\left( 1+x \right)}^{n\beta }},x \right)}^{T}}{{\tilde{D}}^{-1}}{{\tilde{I}}_{n}}. \qquad\square $$

For example, 
$${{A}_{4}}=\left( \begin{matrix}
   1 & \text{  }0 & 0 & 0  \\
   -3 & 1 & 0 & 0  \\
   3 & -2 & 1 & 0  \\
   -1 & 1 & -1\text{  } & 1  \\
\end{matrix} \right)\left( \begin{matrix}
   1 & 0 & 0 & 0  \\
   0 & 2 & 0 & 0  \\
   0 & 0 & 3 & 0  \\
   0 & 0 & 0 & 4  \\
\end{matrix} \right)\left( \begin{matrix}
   1 & 4 & 6 & 4  \\
   0 & 1 & 4 & 6  \\
   0 & 0 & 1 & 4  \\
   0 & 0 & 0 & 1  \\
\end{matrix} \right)\left( \begin{matrix}
   1 & 0 & 0 & 0  \\
   0 & \frac{1}{2} & 0 & 0  \\
   0 & 0 & \frac{1}{3} & 0  \\
   0 & 0 & 0 & \frac{1}{4}  \\
\end{matrix} \right)\left( \begin{matrix}
   1 & 0 & 0 & 0  \\
   3 & 1 & 0 & 0  \\
   3 & 2 & 1 & 0  \\
   1 & 1 & 1 & 1  \\
\end{matrix} \right).$$
{\bfseries Theorem 9.3.}
$$T_{n}^{\beta }{{x}^{p}}=\sum\limits_{m=p}^{n-1}{\frac{1}{\left( m+1 \right)}\left( \begin{matrix}
   n-1-p  \\
   n-1-m  \\
\end{matrix} \right)}{{\left( 1-x \right)}^{n-m-1}}{{t}_{m}}\left( -n\beta +m+1|n\beta ,x \right).$$
{\bfseries Proof.} In this case $p=0$, $1$, … , $n-1$. Then
$$\left[ {{x}^{m}} \right]\tilde{V}_{p+1}^{-1}\tilde{D}{{\left( {{\left( 1+x \right)}^{n\beta }},x \right)}^{T}}{{x}^{p}}=$$
$$=\sum\limits_{i=0}^{m}{{{\left( -1 \right)}^{m-i}}}\left( \begin{matrix}
   p-i  \\
   m-i  \\
\end{matrix} \right)\left( \begin{matrix}
   n\beta   \\
   p-i  \\
\end{matrix} \right)\left( i+1 \right)\frac{\left( n\beta +m-p \right)!}{\left( n\beta +m-p \right)!}=$$
$$=\left( \begin{matrix}
   n\beta   \\
   p-m  \\
\end{matrix} \right)\sum\limits_{i=0}^{m}{{{\left( -1 \right)}^{m-i}}\left( i+1 \right)}\left( \begin{matrix}
   n\beta +m-p  \\
   m-i  \\
\end{matrix} \right)=$$
$$=\left( \begin{matrix}
   n\beta   \\
   p-m  \\
\end{matrix} \right){{\left( -1 \right)}^{m}}\left( \begin{matrix}
   n\beta +m-p-2  \\
   m  \\
\end{matrix} \right)=\left( \begin{matrix}
   n\beta   \\
   p-m  \\
\end{matrix} \right)\left( \begin{matrix}
   -n\beta +p+1  \\
   m  \\
\end{matrix} \right),$$
$$\tilde{V}_{n}^{-1}\tilde{D}{{\left( {{\left( 1+x \right)}^{n\beta }},x \right)}^{T}}{{x}^{p}}={{\left( 1-x \right)}^{n-p-1}}{{t}_{p}}\left( -n\beta +p+1|n\beta ,x \right).$$
$${{\tilde{D}}^{-1}}{{\tilde{V}}_{n}}{{x}^{p}}=\sum\limits_{m=p}^{n-1}{\frac{1}{\left( m+1 \right)}\left( \begin{matrix}
   n-1-p  \\
   n-1-m  \\
\end{matrix} \right)}{{x}^{m}}. \qquad\square $$

In particular, 
$$A_{n}^{\beta }{{x}^{n-1}}=\frac{1}{n}\sum\limits_{m=0}^{n-1}{\left( \begin{matrix}
   n\left( 1-\beta  \right)  \\
   m  \\
\end{matrix} \right)\left( \begin{matrix}
   n\beta   \\
   n-1-m  \\
\end{matrix} \right){{x}^{m}}}, \quad A_{n}^{\beta }{{x}^{0}}=\frac{1}{n}\sum\limits_{m=0}^{n-1}{\left( \begin{matrix}
   -n\beta   \\
   m  \\
\end{matrix} \right)\left( \begin{matrix}
   n\left( 1+\beta  \right)  \\
   n-1-m  \\
\end{matrix} \right){{x}^{m}}},$$
that corresponds to the formula (2).

Denote $\left[ {{P}_{n}}\left( x \right) \right]{{\left( 1,x{}_{\left( \beta  \right)}a\left( x \right) \right)}_{E}}{{=}_{\left( \beta  \right)}}{{\varphi }_{n}}\left( x \right)$.We introduce the matrices $T_{n}^{\beta }={{\tilde{F}}_{n}}{{E}^{n\beta }}\tilde{F}_{n}^{-1}$. Then $T_{n}^{\beta }{{\tilde{\varphi }}_{n}}\left( x \right)={}_{\left( \beta  \right)}{{\tilde{\varphi }}_{n}}\left( x \right)$. For example,
$${{T}_{2}}=\frac{1}{2}\left( \begin{matrix}
   3 & 1  \\
   -1 & 1  \\
\end{matrix} \right),  \quad{{T}_{3}}=\frac{1}{5}\left( \begin{matrix}
   12 & 4 & 1  \\
   -9 & 2 & 3  \\
   2 & -1 & 1  \\
\end{matrix} \right),  \quad{{T}_{4}}=\frac{1}{14}\left( \begin{matrix}
   55 & {55}/{3}\; & 5 & 1  \\
   -66 & 0 & 10 & 6  \\
   30 & -6 & 0 & 6  \\
   -5 & {5}/{3}\; & -1 & 1  \\
\end{matrix} \right).$$
{\bfseries Theorem 9.4.} 
$$T_{n}^{-\beta }={{\tilde{J}}_{n}}T_{n}^{\beta }{{\tilde{J}}_{n}}.$$
{\bfseries Proof.} By the Theorem 8.3.
$${{\tilde{J}}_{n}}{{\tilde{F}}_{n}}{{E}^{n\beta }}\tilde{F}_{n}^{-1}{{\tilde{J}}_{n}}={{\tilde{F}}_{n}}{{E}^{n}}\left( 1,-x \right){{E}^{n\beta }}{{E}^{n}}\left( 1,-x \right)\tilde{F}_{n}^{-1}={{\tilde{F}}_{n}}{{E}^{-n\beta }}\tilde{F}_{n}^{-1}. \qquad\square $$

Matrix $T_{n}^{\beta }$ can be represented as
$$T_{n}^{\beta }={{\tilde{S}}_{n}}A_{n}^{\beta }\tilde{S}_{n}^{-1}=\tilde{V}_{n}^{-1}{{\tilde{C}}_{n}}\tilde{D}{{\left( {{\left( 1+x \right)}^{n\beta }},x \right)}^{T}}{{\tilde{D}}^{-1}}\tilde{C}_{n}^{-1}{{\tilde{V}}_{n}},$$
$${{\tilde{C}}_{n}}\tilde{D}{{x}^{p}}=\left( n+1 \right)!\left( \begin{matrix}
   n+1+p  \\
   p  \\
\end{matrix} \right){{x}^{p}}.$$
For example, 
$${{T}_{4}}=\left( \begin{matrix}
   1 & 0 & 0 & 0  \\
   -3 & 1 & 0 & 0  \\
   3 & -2 & 1 & 0  \\
   -1 & 1 & -1 & 1  \\
\end{matrix} \right)\left( \begin{matrix}
   1 & 0 & 0 & 0  \\
   0 & 6 & 0 & 0  \\
   0 & 0 & 21 & 0  \\
   0 & 0 & 0 & 56  \\
\end{matrix} \right)\left( \begin{matrix}
   1 & 4 & 6 & 4  \\
   0 & 1 & 4 & 6  \\
   0 & 0 & 1 & 4  \\
   0 & 0 & 0 & 1  \\
\end{matrix} \right)\left( \begin{matrix}
   1 & 0 & 0 & 0  \\
   0 & \frac{1}{6} & 0 & 0  \\
   0 & 0 & \frac{1}{21} & 0  \\
   0 & 0 & 0 & \frac{1}{56}  \\
\end{matrix} \right)\left( \begin{matrix}
   1 & 0 & 0 & 0  \\
   3 & 1 & 0 & 0  \\
   3 & 2 & 1 & 0  \\
   1 & 1 & 1 & 1  \\
\end{matrix} \right).$$
{\bfseries Theorem 9.5.}
$$T_{n}^{\beta }{{x}^{p}}=\sum\limits_{m=p}^{n-1}{{{\left( \begin{matrix}
   n+1+m  \\
   m  \\
\end{matrix} \right)}^{-1}}\left( \begin{matrix}
   n-1-p  \\
   n-1-m  \\
\end{matrix} \right)}{{\left( 1-x \right)}^{n-m-1}}{{t}_{m}}\left( -n\beta +n+m+1|n\beta ,x \right).$$
{\bfseries Proof.}
$$\frac{1}{\left( n+1 \right)!}\left[ {{x}^{m}} \right]\tilde{V}_{p+1}^{-1}{{\tilde{C}}_{n}}\tilde{D}{{\left( {{\left( 1+x \right)}^{n\beta }},x \right)}^{T}}{{x}^{p}}=$$
$$=\sum\limits_{i=0}^{m}{{{\left( -1 \right)}^{m-i}}}\left( \begin{matrix}
   p-i  \\
   m-i  \\
\end{matrix} \right)\left( \begin{matrix}
   n\beta   \\
   p-i  \\
\end{matrix} \right)\left( \begin{matrix}
   n+1+i  \\
   i  \\
\end{matrix} \right)\frac{\left( n\beta +m-p \right)!}{\left( n\beta +m-p \right)!}=$$
$$=\left( \begin{matrix}
   n\beta   \\
   p-m  \\
\end{matrix} \right)\sum\limits_{i=0}^{m}{{{\left( -1 \right)}^{m-i}}\left( \begin{matrix}
   n+1+i  \\
   i  \\
\end{matrix} \right)}\left( \begin{matrix}
   n\beta +m-p  \\
   m-i  \\
\end{matrix} \right)=$$
$$=\left( \begin{matrix}
   n\beta   \\
   p-m  \\
\end{matrix} \right){{\left( -1 \right)}^{m}}\left( \begin{matrix}
   n\beta +m-p-n-2  \\
   m  \\
\end{matrix} \right)=\left( \begin{matrix}
   n\beta   \\
   p-m  \\
\end{matrix} \right)\left( \begin{matrix}
   -n\beta +n+p+1  \\
   m  \\
\end{matrix} \right),$$
$$\tilde{V}_{n}^{-1}{{\tilde{C}}_{n}}\tilde{D}{{\left( {{\left( 1+x \right)}^{n\beta }},x \right)}^{T}}{{x}^{p}}=\left( n+1 \right)!{{\left( 1-x \right)}^{n-p-1}}{{t}_{p}}\left( -n\beta +n+p+1|n\beta ,x \right),$$
$${{\tilde{D}}^{-1}}\tilde{C}_{n}^{-1}{{\tilde{V}}_{n}}{{x}^{p}}=\frac{1}{\left( n+1 \right)!}\sum\limits_{m=p}^{n-1}{{{\left( \begin{matrix}
   n+1+m  \\
   m  \\
\end{matrix} \right)}^{-1}}\left( \begin{matrix}
   n-1-p  \\
   n-1-m  \\
\end{matrix} \right)}{{x}^{m}}. \qquad\square $$

In particular, 
$$T_{n}^{\beta }{{x}^{n-1}}={{\left( \begin{matrix}
   2n  \\
   n-1  \\
\end{matrix} \right)}^{-1}}\sum\limits_{m=0}^{n-1}{\left( \begin{matrix}
   n\left( 2-\beta  \right)  \\
   m  \\
\end{matrix} \right)\left( \begin{matrix}
   n\beta   \\
   n-1-m  \\
\end{matrix} \right){{x}^{m}}},$$
$$T_{n}^{\beta }{{x}^{0}}={{\left( \begin{matrix}
   2n  \\
   n-1  \\
\end{matrix} \right)}^{-1}}\sum\limits_{m=0}^{n-1}{\left( \begin{matrix}
   -n\beta   \\
   m  \\
\end{matrix} \right)\left( \begin{matrix}
   n\left( 2+\beta  \right)  \\
   n-1-m  \\
\end{matrix} \right){{x}^{m}}}.$$
that corresponds to the formula (3).

E-mail: {evgeniy\symbol{"5F}burlachenko@list.ru}
\end{document}